\documentclass[review]{elsarticle}

\usepackage{lineno,hyperref}
\modulolinenumbers[5]

\journal{arXiv.org}
\usepackage[margin=1in]{geometry}

\usepackage{latexsym}

\graphicspath{Figures/}
\usepackage[cmex10]{amsmath}

\usepackage{amsfonts}
\usepackage{amsthm}
\usepackage{algorithm}
\usepackage{algpseudocode}
\usepackage{array}
\usepackage{array}

\usepackage{bm}   
\usepackage{multirow} 
\usepackage{bigdelim}
 
\usepackage{booktabs}
\def\NoNumber#1{{\def\alglinenumber##1{}\State #1}\addtocounter{ALG@line}{-1}}
\usepackage{mathtools}
\DeclarePairedDelimiter\ceil{\lceil}{\rceil}

\newtheorem{theorem}{Theorem}

\newtheorem{lemma}{Lemma}
\newtheorem{assumption}{Assumption}
\newtheorem{nameddefinition}{Definition}
\usepackage{xpatch}
\makeatletter
\xpatchcmd{\@thm}{\thm@headpunct{.}}{\thm@headpunct{}}{}{}
\makeatother

  \usepackage[caption=false,font=normalsize,labelfont=sf,textfont=sf]{subfig}

\usepackage{stfloats}

\begin{document}

\begin{frontmatter}

\title{Asynchronous ADMM for Distributed Non-Convex Optimization in Power Systems}

\author[mymainaddress]{Junyao Guo \corref{mycorrespondingauthor}}
\cortext[mycorrespondingauthor]{Corresponding author}
\ead[url]{junyaog@andrew.cmu.edu}

\author[mysecondaryaddress]{Gabriela Hug}
\ead[url]{ghug@ethz.ch}

\author[mymainaddress]{Ozan Tonguz }
\ead[url]{tonguz@ece.cmu.edu}

\address[mymainaddress]{Department of Electrical and Computer Engineering, Carnegie Mellon University, Pittsburgh, PA, USA}
\address[mysecondaryaddress]{Power Systems Laboratory, ETH Z{\"u}rich, Z{\"u}rich, Switzerland }

\begin{abstract}
Large scale, non-convex optimization problems arising in many complex networks such as the power system call for efficient and scalable distributed optimization algorithms. Existing distributed methods are usually iterative and require synchronization of all workers at each iteration, which is hard to scale and could result in the under-utilization of computation resources due to the heterogeneity of the subproblems. To address those limitations of synchronous schemes, this paper proposes an asynchronous distributed optimization method based on the Alternating Direction Method of Multipliers (ADMM) for non-convex optimization. The proposed method only requires local communications and allows each worker to perform local updates with information from a subset of but not all neighbors. We provide sufficient conditions on the problem formulation, the choice of algorithm parameter and network delay, and show that under those mild conditions, the proposed asynchronous ADMM method asymptotically converges to the KKT point of the non-convex problem. We validate the effectiveness of asynchronous ADMM by applying it to the Optimal Power Flow problem in multiple power systems and show that the convergence of the proposed asynchronous scheme could be faster than its synchronous counterpart in large-scale applications.
\end{abstract}


\end{frontmatter}


\section{Introduction}
\label{intro}
The future power system is expected to integrate large volumes of distributed generation resources, distributed storages, sensors and measurement units, and flexible loads. As the centralized architecture will be prohibitive for collecting measurements from all of those newly installed devices and coordinating them in real-time, it is expected that the grid will undergo a transition towards a more distributed control architecture. Such transition calls for distributed algorithms for optimizing the management and operation of the grid that does not require centralized control and therefore scales well to large networks. 

Major functions for power system management include Optimal Power Flow (OPF), Economic Dispatch, and State Estimation, which can all be abstracted as the following optimization problem in terms of variables assigned to different control regions \cite{kar2014distributed}:
\begin{subequations}
\label{eq:intro}
\begin{align}
\underset{{x}} {\text{minimize}} ~~~&\sum_{k} f_{k}(x_{k})\\
\label{eq:noncouplingconstraint}
\text{subject to}~~~& n({x}_{k})\leq 0,~\forall k\\
\label{eq:couplingconstraint}
& c(x_{1},..., x_{k},..., x_{K}) \leq 0,
\end{align}
\end{subequations}
\normalsize
where $K$ and $x_k$ denote the total number of regions and the variables in region $k$, respectively. Functions $f_{k}(\cdot)$, $n(\cdot)$ and $c(\cdot)$ are smooth but possibly non-convex functions. Variable $x$ is bounded due to the operating limits of the devices; that is, $x \in \mathcal{X}$ where $\mathcal{X}$ is a compact smooth manifold. Note that constraint (\ref{eq:couplingconstraint}) is usually referred to as the coupling constraint as it includes variables from multiple control regions and therefore couples the updates of those regions. 

In this paper, we are interested in developing efficient optimization algorithms to solve problem (\ref{eq:intro}) in a distributed manner. Many iterative distributed optimization algorithms have been proposed for parallel and distributed computing \cite{conejo2006decomposition}\cite{bertsekas1989parallel} where among these algorithms, it has been shown that the Alternating Direction Method of Multipliers (ADMM) often exhibits good performance for non-convex optimization. 
The convergence of ADMM in solving non-convex problems is characterized in two recent studies \cite{magnusson2016convergence}\cite{wang2015global} with a synchronous implementation where all subproblems need to be solved before a new iteration starts. In fact, the majority of distributed algorithms are developed based on the premise that the workers that solve subproblems are synchronized. However, synchronization may not be easily obtained in a distributed system without centralized coordination, which to a certain degree defeats the purpose of the distributed algorithm. Moreover, the sizes and complexities of subproblems are usually dependent on the system's physical configuration, and therefore are heterogeneous and require different amounts of computation time. Therefore, even if synchronization is achievable, it may not be the most efficient way to implement distributed algorithms. Furthermore, the communication delays among workers are also heterogeneous which are determined by the communication infrastructures and technologies used. In a synchronous setting, all workers need to wait for the slowest worker to finish its computation or communications. This may lead to the under-utilization of both the computation and communication resources as some workers remain idle for most of the time.

To alleviate the above limitations of synchronous schemes, in this paper, we propose a distributed asynchronous optimization approach which is based on the state-of-the-art ADMM method \cite{boyd2011distributed}. We extend this method to fit into an asynchronous framework where a message-passing model is used and each worker is allowed to perform local updates with partial but not all updated information received from its neighbors. Particularly, the proposed method is scalable because it only requires local information exchange between neighboring workers but no centralized or master node. The major contribution of this paper is to show that the proposed asynchronous ADMM algorithm asymptotically satisfies the first-order optimality conditions of problem (\ref{eq:intro}), under the assumption of bounded delay of the worker and some other mild conditions on the objective function and constraints. To the best of our knowledge, this is the first time that distributed ADMM is shown to be convergent for a problem with non-convex coupling constraints (see (\ref{eq:intro})) under an asynchronous setting. Also, we show that the proposed asynchronous scheme can be applied to solving the AC Optimal Power Flow (AC OPF) problem in large-scale transmission networks, which provides a promising tool for the management of future power grids that are likely to have a more distributed architecture.

The rest of the paper is organized as follows: Section \ref{related work} presents related work and contrasts our contributions. The standard synchronous ADMM method is first sketched in Section \ref{synchronousADMM}, while the asynchronous ADMM method is given in Section \ref{asynchronousADMM}. The sufficient conditions and the main convergence result are stated in Section \ref{convanalysis} and the proof of convergence is shown in Section \ref{maintheoremproof}. Section \ref{results} demonstrates numerical results on the performance of asynchronous ADMM for solving the AC OPF problem. Finally, Section \ref{conclusion} concludes the paper and proposes future studies.

\section{Related Works}
\label{related work}
The synchronization issue has been systematically studied in the research fields of distributed computing with seminal works \cite{bertsekas1989parallel}\cite{dwork1988consensus}. While the concept of asynchronous computing is not new, it remains an open question whether those methods can be applied to solving non-convex problems. Most of the asynchronous computing methods proposed can only be applied to convex optimization problems \cite{bertsekas1989parallel}\cite{peng2016arock}. These methods therefore can only solve convex approximations of the non-convex problems which may not be exact for all types of systems \cite{7285913, abboud2014asynchronous,nguyen2016distributed,chang2017scheduled}. Recently, there are a few asynchronous algorithms proposed that tackle problems with some level of non-convexity. Asynchronous distributed ADMM approaches are proposed in \cite{7423789} and \cite{kumar2017asynchronous} for solving consensus problems with non-convex local objective functions and convex constraint sets; but, the former approach requires a master node and the latter uses linear approximations for solving local non-convex problems. In \cite{cannelli2016asynchronous}, a probablistic model for a parallel asynchronous computing method is proposed for minimizing non-convex objective functions with separable non-convex sets. However, none of the aforementioned studies handles non-convex coupling constraints that include variables from multiple workers. The non-convex problem (see (\ref{eq:intro})) studied in this work does include such constraints and our approach handles them without convex approximations. 

A further difference concerns the communication graph and the information that is exchanged. A classical problem studied in most research is the consensus problem where the workers are homogeneous and they aim to find a common system parameter. The underlying network topology is either a full mesh network where any pair of nodes can communicate \cite{dwork1988consensus} or a star topology with a centralized coordinator \cite{7423789}\cite{zhang2014asynchronous}. Different from the consensus problem, we consider partition-based optimization where the workers represent regions or subnetworks with different sizes. Furthermore, each worker only communicates with its physically connected neighbors and the information to be exchanged only contains the boundary conditions but not all local information. Thereby, the workers are heterogeneous and the communication topology is a partial mesh network with no centralized/master node needed.

\section{Synchronous Distributed ADMM}
\label{synchronousADMM}
Geographical decomposition of the system is considered in this paper where a network is partitioned into a number of subnetworks each assigned to a worker for solving the local subproblem; i.e., worker $k$ updates a subset $x_k$ of the variables. The connectivity of the workers is based on the network topology; i.e., we say two workers are neighbors if their associated subnetworks are physically connected and some variables of their variable sets appear in the same coupling constraints. In the following analysis, we use $\mathcal{N}_{k}$ and $\mathcal{T}$ to denote the set of neighbors that connect to worker $k$ and the set of edges between any pair of neighboring workers. 

To apply the distributed ADMM approach to problem (\ref{eq:intro}), we introduce auxiliary variables $z_k = \{z_{k,l} | \forall l \in \mathcal{N}_k\}$ for each worker $k$ to denote the boundary conditions that neighboring workers should agree on. Then problem (\ref{eq:intro}) can be expressed as follows:
\begin{subequations}
\label{eq:OPF}
\begin{align}
\label{diseq1}
\underset{{x}, {z}} {\text{minimize}} ~~~&\sum_{k} F_{k}(x_{k})\\
\label{diseq2}
\text{subject to}~~~& A_{k}{x}_{k}={z}_{k},~\forall k\\
\label{diseq4}
& z_{k,l} = z_{l,k}, ~\forall (k,l) \in \mathcal{T},
\end{align}
\end{subequations}
\normalsize
where $F_{k}(x_{k}) =  f_{k}({x}_{k})+\eta_{\mathcal{X}_k}(x_{k}) $ denotes the local objective function and $\eta_{\mathcal{X}}(\cdot)$ is the indicator function of set $\mathcal{X}$, with $\eta_{\mathcal{X}}(x) = 0$ if $x \in {\mathcal{X}}$ and $+\infty$ if $x \notin {\mathcal{X}}$. Constraint (\ref{diseq2}) establishes the relations between $x_k$ and $z_k$ and constraint (\ref{diseq4}) enforces the agreement on the boundary conditions of neighboring workers. Note that by choosing $A_k$ as the identity matrix, problem (\ref{eq:OPF}) reduces to the standard consensus problem where all workers should find a common variable $z$. Here we allow $A_k$ to not have full column rank; i.e., the neighboring workers only need to find common values for the boundary variables but do not need to share information on all local variables, which greatly reduces the amount of information to be exchanged among workers. 

The ADMM algorithm minimizes the Augmented Lagrangian function of (\ref{eq:OPF}), which is given as follows \cite{boyd2011distributed}:
\begin{equation}
\begin{aligned}
L({x},{z},{\lambda})=&\sum_{k} \Big\{ F_{k}(x_{k})+{\lambda}_{k}^{\top}(A_{k}{x}_{k}-{z}_{k})\\
&~~~+\frac{\rho}{2}\|A_{k}{x}_{k}-{z}_{k}\|^{2}\Big\} + \eta_{\mathcal{Z}}(z),
\end{aligned}
\label{eq:Aug}
\end{equation}
\normalsize
where $z$ denotes the superset of all auxiliary variables $z_k, \forall k$ and $\mathcal{Z}$ denotes the feasible region of $z$ imposed by constraint (\ref{diseq4}). The standard synchronous ADMM method minimizes (\ref{eq:Aug}) by iteratively carrying out the following updating steps \cite{boyd2011distributed}:
\begin{subequations}
\label{eq:ADMMIter}
\begin{align}
\label{eq:ADMMIterzupdate}
z-\text{update}: & ~~~~{z}^{\nu+1}=\text{argmin}~~L({x}^{\nu},{z},{\lambda}^{\nu})\\
\label{eq:ADMMIterxupdate}
x-\text{update}: &~~~~{x}^{\nu+1} = \text{argmin}~~ L({x},{z}^{\nu+1},{\lambda}^{\nu})\\
\label{eq:ADMMIterlambdaupdate}
\lambda-\text{update}: &~~~~{\lambda}^{\nu+1}={\lambda}^{\nu}+{\rho}(A{x}^{\nu+1}-{z}^{\nu+1}),
\end{align}
\normalsize
\end{subequations}
where $\nu$ denotes the counter of iterations. With ${z}$ fixed, each subproblem in the $x$-update only contains the local variables ${x}_{k}$, such that the subproblems can be solved independently of each other. The $\lambda$-update can also be performed locally. The $z$-update requires the information from two neighboring workers, thus can also be carried out locally as long as the information from neighboring workers is received.



We define the residue of ADMM as
\begin{equation}
\Gamma_k^{\nu+1}= \left \|
\begin{aligned} 
A_{k}{x}_{k}^{\nu+1}&-{z}_{k}^{\nu+1} \\
{z}_{k}^{\nu+1}&-{z}_{k}^{\nu}
\end{aligned}
\right \|_{\infty},
\label{eq:residue}
\end{equation}
\normalsize
where the two terms denote the primal residue and dual residue\cite{boyd2011distributed}, respectively. The stopping criterion is defined as that both $\Gamma_k$ and the maximum constraint mismatch for all workers are smaller than some $\epsilon$ \cite{boyd2011distributed}. Under the non-convex setting, the convergence of synchronous ADMM to a KKT stationary point $\{x^{\star}, z^{\star}, \lambda^{\star}\}$ is proved in \cite{erseghe2015distributed} with the assumption that both $x$ and $\lambda$ are bounded and that a local minimum can be identified when solving the local subproblems. 
\section{Asynchronous Distributed ADMM}
\label{asynchronousADMM}
\begin{figure}[t]
\setlength{\abovecaptionskip}{0.2cm} 
\centering
\captionsetup[subfigure]{captionskip= 0 cm}
\subfloat[synchronous, $p=1$]
{
\label{fig:illusync}
\includegraphics[trim = 0mm 0mm 0mm 0mm, clip=true,width=7.5cm]{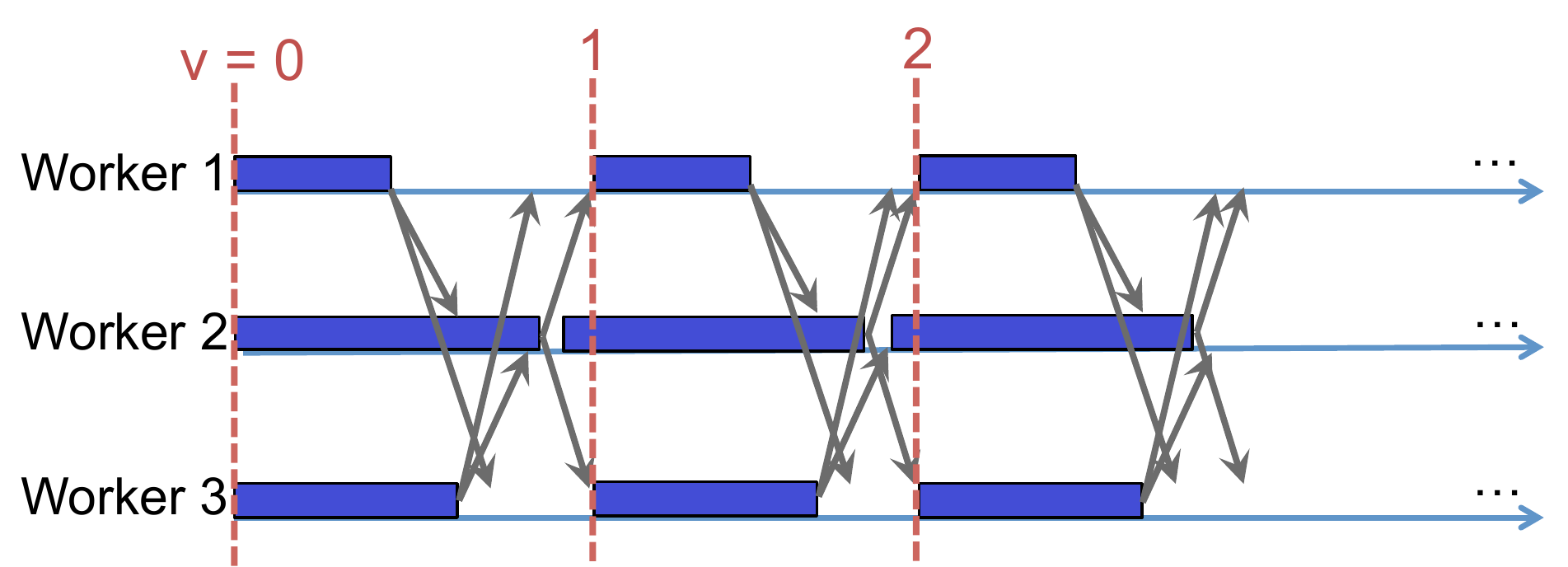} 
}
\hspace{0cm}
\subfloat[asynchronous, $p=0.1$]
{
\label{fig:illuasync}
\includegraphics[trim =0mm 0mm 0mm 0mm, clip=true,width=7.5cm]{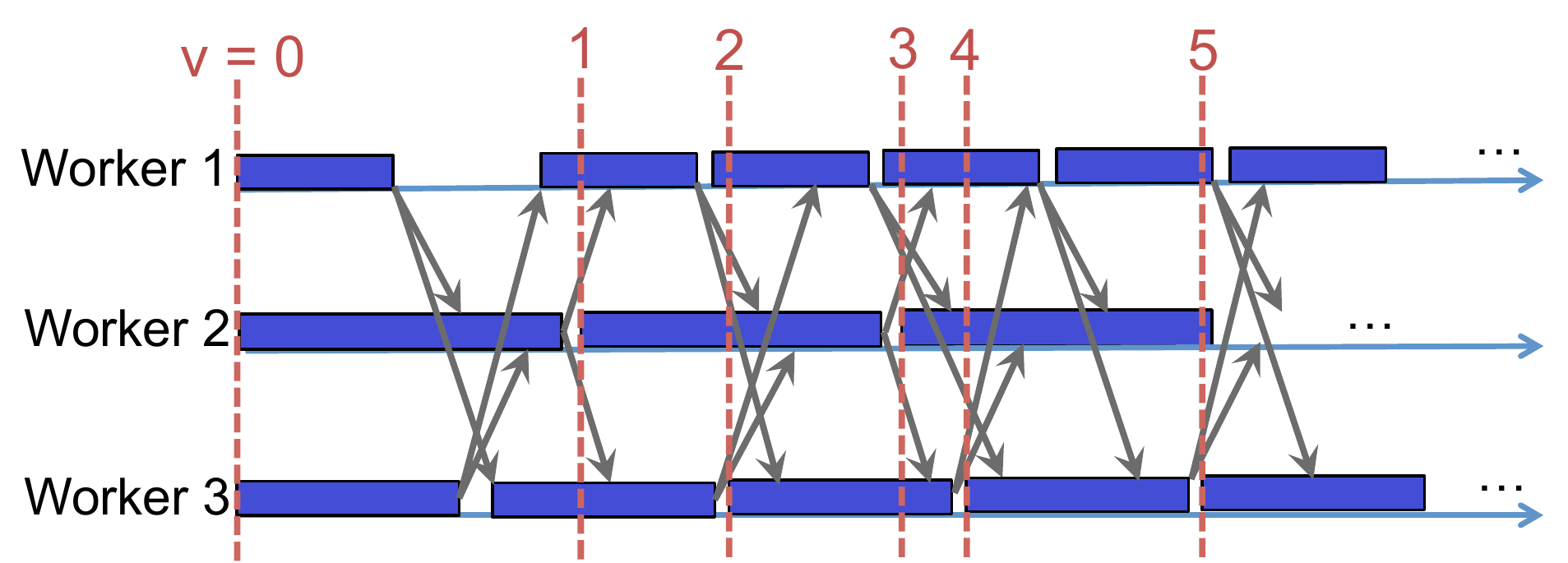} 
}
\caption{Illustration of synchronous and asynchronous distributed ADMM.}
\label{fig:asyncillustration}
\end{figure}

Now, we extend the synchronous ADMM into an asynchronous version where each worker determines when to perform local updates based on the messages it receives from neighbors. We say that a neighbor $l$ has `arrived' at worker $k$ if the updated information of $l$ is received by $k$. We assume partial asynchrony \cite{bertsekas1989parallel} where the delayed number of iterations of each worker is bounded by a finite number. We introduce a parameter $p$ with $0<p\leq1$ to control the level of asynchrony. Worker $k$ will update its local variables after it receives new information from at least $\ceil{p|\mathcal{N}_{k}|}$ neighbors with $|\mathcal{N}_{k}|$ denoting the number of neighbors of worker $k$. In the worst case, any worker should wait for at least one neighbor because otherwise its local update will make no progress as it has no new information. Figure \ref{fig:asyncillustration} illustrates the proposed asynchronous scheme by assuming three workers, each connecting to the other two workers. The blue bar denotes the local computation and the grey line denotes the information passing among workers. The red dotted line marks the count of iterations, which will be explained in Section \ref{convanalysis}. As shown in Fig. \ref{fig:illusync}, synchronous ADMM can be implemented by setting $p = 1$, where each worker does not perform local computation until all neighbors arrive. Figure \ref{fig:illuasync} shows an asynchronous case where each worker can perform its local update with only one neighbor arrived, which could reduce the waiting time for fast workers. In the rest of the paper, we do not specify the setting of $p$ but only consider it to be a very small value such that each worker can perform its local update as long as it received information from at least one neighbor, which indeed represents the highest level of asynchrony.

Algorithm \ref{algADMMasync} presents the asynchronous ADMM approach from each region's perspective with $\nu_{k}$ denoting the local iteration counter. Notice that for the $z$-update, we add a proximal term $\frac{\alpha}{2}\|z_{k} -z_{k}^{\nu_{k}}\|^{2}$ with $\alpha \geq 0$ which is a sufficient condition for proving the convergence of ADMM under asynchrony. The intuition of adding this proximal term is to reduce the stepsize of updating variables to offset the error brought by asynchrony. Also, in the $z$-update, only the entries corresponding to arrived neighbors in $z_{k}$ are updated.
\begin{algorithm}[t]
\caption{Asynchronous ADMM in region $k$}
\label{algADMMasync}
\begin{algorithmic}[1]
\State \textbf{Initialization} 
\NoNumber{Given ${x}_{k}^{0}$, set ${\lambda}_{k}^{0}={0}$, ${\rho}=\rho_{0} $, $\nu_{k}=0, z_{k}^{0} \in \mathcal{Z}$ }
\State \textbf{Repeat}
\State \hspace{\algorithmicindent}\textbf{Update} $x$ by
\begin{equation}
{x}_{k}^{\nu_{k}+1} =  \underset{x_{k}}{\text{argmin}} ~L({x_{k}}, z^{\nu_{k}}_{k},{\lambda}^{\nu_{k}}_{k})\label{eq:xupdatelocal}
\end{equation}
\State \hspace{\algorithmicindent}\textbf{Update}  ${\lambda}$ using
\begin{equation}
{\lambda}_{k}^{\nu_{k}+1}={\lambda}_{k}^{\nu_{k}}+{\rho}(A_{k}{x}_{k}^{\nu_{k}+1}-{z}_{k}^{\nu_{k}})
\end{equation}
\State\hspace{\algorithmicindent}Send $\{A_{k,l}x_{k}^{\nu_{k}+1}, \lambda_{k,l}^{\nu_{k}+1}\}$ to region $l$, $\forall l \in \mathcal{N}_{k}$
\State\hspace{\algorithmicindent}\textbf{Repeat}
\State\hspace{\algorithmicindent}\hspace{\algorithmicindent}\textbf{Wait} until at least $\ceil{p|\mathcal{N}_{k}|}$ neighbors arrive
\State \hspace{\algorithmicindent}\textbf{Update} $z_{k}$ associated with arrived neighbors 
\begin{equation}
z_k^{\nu_{k}+1} =  \underset{z_k}{\text{argmin}}~~L({x}^{\nu_{k}}_{k},z_{k},{\lambda}^{\nu_{k}}_{k})+\frac{\alpha}{2}\|z_{k} -z_{k}^{\nu_{k}}\|^{2}
\label{eq:zupdatelocal}
\end{equation}
\State \hspace{\algorithmicindent}Set $\nu_{k} \leftarrow \nu_{k}+1$
\State \textbf{Until} a predefined stopping criterion is satisfied
\normalsize
\end{algorithmic}
\end{algorithm}

\section{Convergence Analysis}
\label{convanalysis}
For analyzing the convergence property of Algorithm \ref{algADMMasync}, we introduce a global iteration counter $\nu$ and present Algorithm \ref{algADMMasync} from a global point of view in  Algorithm \ref{algADMMasyncglobal} \emph{as if} there is a master node monitoring all the local updates. Note that such master node is not needed for the implementation of Algorithm \ref{algADMMasync} and the global counter only serves the purpose of slicing the execution time of Algorithm \ref{algADMMasync} for analyzing the changes in variables during each time slot. We use the following rules to set the counter $\nu$: 1) $\nu$ can be increased by 1 at time $t_{\nu + 1}$ when some worker is ready to start its local $x$-update; 2) the time period $(t_{\nu}, t_{\nu + 1}]$ should be as long as possible; 3) there is no worker that finishes $x$-update more than once in $(t_{\nu}, t_{\nu + 1}]$; 4) $\nu$ should be increased by 1 before any worker receives new information after it has started its $x$-update. The third rule ensures that each $x$-update is captured in one individual iteration and the fourth rule ensures that the $z$ used for any $x$-update during $(t_{\nu}, t_{\nu + 1}]$ is equal to the $z$ measured at $t_{\nu + 1}$. This global iteration counter is represented by red dotted lines in Fig. \ref{fig:asyncillustration}. 

We define $\mathcal{A}_{\nu}$ as the index subset of workers who finishes $x$-updates during the time $(t_{\nu}, t_{\nu + 1}]$ with $0 \leq |\mathcal{A}_{\nu}| \leq K$. Note that with $|\mathcal{A}_{\nu}| = K$, Algorithm \ref{algADMMasyncglobal} is equivalent to synchronous ADMM. We use $\mathcal{T}_{(t_{\nu - 1}, t_{\nu}]}$ to denote the set of workers that exchange information at iteration $\nu$; i.e., $(k,l) \in \mathcal{T}_{(t_{\nu - 1}, t_{\nu}]}$ denotes that the updated information from worker $k$ arrives at $l$ during the time $(t_{\nu - 1}, t_{\nu}]$. 
\begin{algorithm}[t]
\caption{Asynchronous ADMM from a global view}
\label{algADMMasyncglobal}
\begin{algorithmic}[1]
\State \textbf{Initialization} 
\NoNumber{Given ${x}^{0}$, $\rho$, set ${\lambda}^{0}={0}$, $\nu=0, z^{0} \in \mathcal{Z}$ }
\State \textbf{Repeat}
\State \hspace{\algorithmicindent}\textbf{Update} 
\NoNumber{
\begin{equation}
{x}_{k}^{\nu+1} = \left \{
\begin{aligned} 
&\text{argmin} ~F_{k}({x}_{k})+{\lambda}_{k}^{\nu\top}A_{k}{x}_{k}+\frac{{\rho}}{2}\|A_{k}{x}_{k}-{z}_{k}^{\bar{\nu}_{k}+1}\|^{2} & \text{if}~k \in \mathcal{A}_{\nu}  \\
&x_{k}^{\nu}  & \text{otherwise}
\end{aligned}
\right.
\label{eq:xupdateglobal}
\end{equation}}
\NoNumber{
\begin{equation}
 {\lambda}_{k}^{\nu+1} = \left \{
\begin{aligned} 
&{\lambda}_{k}^{\nu}+ {\rho}(A_{k}{x}_{k}^{\nu+1}-{z}_{k}^{\bar{\nu}_{k}+1})& \text{if}~k \in \mathcal{A}_{\nu}  \\
&{\lambda}_{k}^{\nu}  & \text{otherwise}
\end{aligned}
\right.
\label{eq:lambdaupdateglobal}
\end{equation}}
\NoNumber{For $\forall k \in \mathcal{A}_{\nu}, l \in \mathcal{N}_{k}~~ \text{and}   ~~\forall (l,k) \in \mathcal{T}_{(t_{\nu - 1}, t_{\nu}]}$:}
\NoNumber{
\begin{equation}
\begin{aligned}
z_{k,l}^{\nu + 1} = &\bigg(\lambda_{k,l}^{\nu + 1}+\lambda_{l,k}^{\nu + 1}+\rho A_{k,l}x_{k}^{\nu + 1} +\rho A_{l,k}x_{l}^{\nu + 1}+\alpha z_{k,l}^{\nu } \bigg) / (2\rho+\alpha)
\end{aligned}
\label{eq:zupdateglobal}
\end{equation}}
\State \hspace{\algorithmicindent}Set $\nu \leftarrow \nu +1$
\State \textbf{Until} a predefined stopping criterion is satisfied
\normalsize
\end{algorithmic}
\end{algorithm}
Now, we formally introduce the assumption of partial asynchrony (bounded delay).
\begin{assumption}
Let $\omega > 0$ be a maximum number of global iterations between the two consecutive $x$-updates for any worker $k$; i.e., for all $k$ and global iteration $\nu > 0$, it must hold that $k\in \mathcal{A}_{\nu} \cup \mathcal{A}_{\nu-1} \cdots \cup \mathcal{A}_{\max \{\nu - \omega + 1, 0\}}$
\label{assumption:boundeddelay}
\end{assumption}
Define $\bar{\nu}_{k}$ as the iteration number at the start of the $x$-update that finishes at iteration $\nu$. Then, under Assumption \ref{assumption:boundeddelay} and due to the fact that any worker can only start a new $x$-update after it has finished its last $x$-update, it must hold that
\begin{equation}
\max\{\nu - \omega, 0\} \leq \bar{\nu}_{k} < \nu, ~~~\forall \nu > 0.
\label{eq:iterdelaybound}
\end{equation}
\normalsize

The $z$-update (\ref{eq:zupdateglobal}) is derived from the optimality condition of (\ref{eq:zupdatelocal}). As an example, we show how to update variable $z_{k,l}$ (same for $z_{l,k}$) for one pair of neighboring workers $k$ and $l$. To fulfill $z_{k,l} = z_{l,k}$, we substitute $z_{l,k}$ with $z_{k,l}$ and then remove the part $\eta_{\mathcal{Z}}(z)$ from (\ref{eq:Aug}). The remaining part that contains $z_{k,l}$ in (\ref{eq:Aug}) can be written as:
\begin{equation}
\begin{aligned}
 &L' (x_{k}^{\nu+1}, z_{k,l},\lambda_{k}^{\nu+1})  = -(\lambda_{k,l}^{\nu+1}+  \lambda_{l,k}^{\nu+1})z_{k,l} 
 + \frac{\rho}{2}\|A_{k,l}x_{k}^{\nu+1}-z_{k,l} \|^{2} +
 \frac{\rho}{2}\|A_{l,k}x_{l}^{\nu + 1}-z_{k,l}\|^{2} +\frac{\alpha}{2}\|z_{k,l}-z_{k,l}^{\nu}\|^{2}
\end{aligned}
\label{eq:Lsingz}
\end{equation}
\normalsize
The optimality condition of (\ref{eq:zupdatelocal}) then yields
\small
\begin{equation}
\begin{aligned}
&\lambda_{k,l}^{\nu + 1} + \lambda_{l,k}^{\nu + 1} + \rho (A_{k,l}x_{k}^{\nu + 1} - z_{k,l}^{\nu+1}) + \rho (A_{l,k}x_{l}^{\nu + 1} - z_{k,l}^{\nu+1}) -\alpha(z_{k,l}^{\nu+1}-z_{k,l}^{\nu}) = 0,
\end{aligned}
\label{eq:zoptimality}
\end{equation}
\normalsize
which results in (\ref{eq:zupdateglobal}). Thereby, worker $k$ will update $z_{k}$ locally once 1) it receives $\lambda_{l}$, $\rho_{l}$ and $A_{l,k}x_{l}$ from $\forall l \in \mathcal{N}_{k}$ or 2) it finishes local $x$ and $\lambda$ updates.

Before we state our main results of convergence analysis, we need to introduce the following definitions and make the following assumptions with respect to problem (\ref{eq:OPF}).
\begin{nameddefinition}[\bf{Restricted prox-regularity}]\cite{wang2015global}
Let $D \in \mathbb{R}_{+}$, $f~:~\mathbb{R}^{N} \rightarrow \mathbb{R} \cup \{\infty\}$, and define the exclusion set
\begin{equation}
S_{D} ~:= \{x \in \text{dom}(f)~:~ \|d\| > D ~~\text{for all} ~~d \in \partial f(x)\}.
\end{equation}
\normalsize
$f$ is called \emph{restricted prox-regular} if, for any $D > 0$ and bounded set $T \subseteq \text{dom}(f)$, there exists $\gamma > 0$ such that 
\begin{equation}
\begin{aligned}
&f(y) + \frac{\gamma}{2}\|x - y\|^{2} \geq f(x) + \langle d, y - x \rangle, \forall x \in T \setminus S_{D}, y \in T, d \in \partial f(x), \|d\| \leq D. 
\end{aligned}
\end{equation}
\normalsize
\label{def:proxregular}
\end{nameddefinition}
\vspace{-0.2cm}
\begin{nameddefinition}[\bf{Strongly convex functions}]
A convex function $f$ is called strongly convex with modulus $\sigma$ if either of the following holds:
\begin{enumerate}
\item there exists a constant $\sigma > 0$ such that the function \small$f(x) - \frac{\sigma}{2}\|x\|^{2}$\normalsize is convex;
\item there exists a constant $\sigma > 0$ such that for any $x, y \in \mathbb{R}^{N}$ we have:
\begin{equation}
f(y) \geq f(x) + \langle f'(x), y - x \rangle + \frac{\sigma}{2}\|y-x\|^{2}.
\end{equation}
\normalsize
\end{enumerate}
\label{def:strongconvexity}
\end{nameddefinition}
\vspace{-0.2cm}
The following assumptions state the desired characteristics of the objective functions and constraints.
\begin{assumption}
$\mathcal{X}$ is a compact smooth manifold and there exists constant $M_{2} > 1$ such that, $\forall \nu_{1}, \nu_{2}$
$$\frac{1}{M_{2}} \|A_{k}x_{k}^{\nu_{1}} -A_{k}x_{k}^{\nu_{2}}\| \leq \|x_{k}^{\nu_{1}} -x_{k}^{\nu_{2}} \|\leq M_{2} \|A_{k}x_{k}^{\nu_{1}} -A_{k}x_{k}^{\nu_{2}} \|.$$\normalsize
\label{assumption:Axmapping}
\end{assumption}
\vspace{-0.2cm}
Assumption \ref{assumption:Axmapping} allows $A_{k}$ to not have full column rank, which is more realistic for region-based optimization applications. Since $\mathcal{X}$ is compact and $\|A_{k}\| < \infty$, Assumption \ref{assumption:Axmapping} is satisfied.

\begin{assumption}
Each function $F_{k}$ is restricted prox-regular (Definition \ref{def:proxregular}) and its subgradient $\partial F_{k}$ is Lipschitz continuous with a Lipschitz constant $M_{1} > 0$.
\label{assumption:proxregular}
\end{assumption}
The objective function $F_{k}$ in problem (\ref{eq:OPF}) includes indicator functions whose boundary is defined by $\mathcal{X}$. Recall that we assume $\mathcal{X}$ is compact and smooth and as stated in \cite{wang2015global}, indicator functions of compact smooth manifolds are restricted prox-regular functions. 

\begin{assumption}
The subproblem (\ref{eq:xupdatelocal}) is feasible and a local minimum $x_{k} \in \mathcal X_{k}$ is found at each $x$-update. 
\label{assumption:xfeasible}
\end{assumption}
\color{black}
Assumption \ref{assumption:xfeasible} can be satisfied if the subproblem is not ill-conditioned and the solver used to solve the subproblem is sufficiently robust to identify a local optimum, which is generally the case observed from our empirical studies.

\begin{assumption}
$A_{k}A_{k}^{\top}$ is invertible for all $k$, and define $B_{k} = (A_{k}A_{k}^{\top})^{-1}A_{k} $. Also, let $\sigma_{\max}(\cdot)$ denote the operator of taking the largest eigenvalue of a symmetric matrix and define $C = \max\{\sigma_{\max}(B_{k}^{\top}B_{k}), \forall k\}$.
\label{assumption:Ainvertible}
\end{assumption}
\begin{assumption}
$\lambda$ is bounded, and \small$$-\infty <L(x^{\nu}, z^{\nu}, \lambda^{\nu}) < \infty, \text{if}~~ x \in \mathcal{X}$$\normalsize
\label{assumption:Lbound}
\end{assumption}
\vspace{-0.4cm}
$\lambda$ can be bounded by the projection onto a compact box, i.e., $\lambda^{\nu} \leftarrow \max(\lambda^{\min}, \min(\lambda^{\nu},\lambda^{\max}))$. Then Assumption \ref{assumption:Lbound} holds as all the terms in $L(x^{\nu}, z^{\nu}, \lambda^{\nu})$ are bounded with $x$ in the compact feasible region.

The main convergence result of asynchronous ADMM is stated below.
\begin{theorem}
Suppose that Assumptions \ref{assumption:boundeddelay} to \ref{assumption:Lbound} hold. Moreover, choose
\begin{equation}
\begin{aligned}
\rho & > (\gamma+CM_{1}^{2})M_{2}^{2}+\sqrt{(\gamma+CM_{1}^{2})^{2}M_{2}^{4}+4CM_{1}^2M_{2}^2},\\
\alpha & > \frac{(2\rho M_{2}^{4}+1)(\omega-1)^{2}}{2}-\rho.
\end{aligned}
\label{eq:rhoalpha}
\end{equation}
\normalsize
Then, $(\{x_{k}^{\nu}\}_{k=1}^{K}, z^{\nu}, \{\lambda_{k}^{\nu}\}_{k=1}^{K})$ generated by (\ref{eq:xupdatelocal}) to (\ref{eq:zupdatelocal}) (or equivalently (\ref{eq:xupdateglobal}) to (\ref{eq:zupdateglobal})) are bounded and have limit points that satisfy the KKT conditions of problem (\ref{eq:OPF}) for local optimality.
\label{thm:ADMMconvergence}
\end{theorem}
\vspace{-0.2cm}
\section{Proof of Theorem \ref{thm:ADMMconvergence}}
\label{maintheoremproof}
The essence of proving Theorem \ref{thm:ADMMconvergence} is to show the sufficient descent of the Augmented Lagrangian function (\ref{eq:Aug}) at each iteration and that the difference of (\ref{eq:Aug}) between two consecutive iterations is summable. The proof of Theorem \ref{thm:ADMMconvergence} uses the following lemmas, which are proved in the Appendix.
\begin{lemma}
Suppose that Assumption \ref{assumption:Axmapping} to \ref{assumption:Ainvertible} hold. Then it holds that
\begin{equation}
\begin{aligned}
&L(x^{\nu+1}, z^{\nu+1}, \lambda^{\nu+1})-L(x^{\nu}, z^{\nu}, \lambda^{\nu})\\
& \leq (\frac{\gamma+CM_{1}^{2}}{2} - \frac{\rho}{4 M_{2}^{2}}+\frac{CM_{1}^{2}}{\rho}) \sum_{k\in\mathcal{A}_{\nu}}\|x_{k}^{\nu+1}-x_{k}^{\nu}\|^2\\&~~~~-(\rho+\alpha)\|z^{\nu + 1} - z^{\nu}\|^{2}+\frac{2\rho M_{2}^{4}+1}{2} \sum_{k\in\mathcal{A}_{\nu}}\|{z}_{k}^{\bar{\nu}_{k}+1} - {z}_{k}^{\nu}\|^{2}.
\end{aligned}
\label{eq:lemmaLdiff}
\end{equation}
\label{lemma:Ldifference}
\end{lemma}
\normalsize
\vspace{-0.2cm}
Due to the term $\|z_{k}^{\bar{\nu}_{k}+1}- z_{k}^{\nu}\|^{2}$ which is caused by the asynchrony of the updates among workers, (\ref{lemma:Ldifference}) is not necessarily decreasing. We bound this term by Lemma \ref{lemma:errorbound}.

\begin{lemma}
Suppose that Assumption \ref{assumption:boundeddelay} holds. Then it holds that 
\begin{equation}
\sum_{\phi = 1}^{\nu}\sum_{k\in\mathcal{A}_{\phi}}\|z_{k}^{\bar{\phi}_{k}+1}-z_{k}^{\phi} \|^{2} \leq 2(\omega-1)^{2} \sum_{\phi = 1}^{\nu} \|z^{\phi + 1}- z^{\phi}\|^{2}.
\label{eq:errorbound}
\end{equation}
\normalsize
\label{lemma:errorbound}
\end{lemma}
\vspace{-0.3cm}
Using Lemma \ref{lemma:Ldifference} and Lemma \ref{lemma:errorbound}, we now prove Theorem \ref{thm:ADMMconvergence}.
\begin{proof}[Proof of Theorem \ref{thm:ADMMconvergence}]
Any KKT point $(\{x_{k}^{\star}\}_{k = 1}^{K}, z^{\star},\{\lambda_{k}^{\star}\}_{k = 1}^{K})$ of problem (\ref{eq:OPF}) should satisfy the following conditions
\begin{subequations}
\label{eq:OPFKKT}
\begin{align}
\label{eq:KKTx}
& \partial F_{k}(x_{k}^{\star}) + A_{k}^{\top}\lambda_{k}^{\star} = \bm{0}, ~~\forall k\\
\label{eq:KKTz}
& \lambda_{k,l}^{\star} + \lambda_{l,k}^{\star} = \bm{0}, ~~\forall (k,l) \in \mathcal{T}\\
\label{eq:KKTlambda}
& A_{k}{x}_{k}^{\star}-{z}_{k}^{\star}=\bm{0},~~\forall k
\end{align}
\label{eq:KKTall}
\end{subequations}
\normalsize
By taking the telescoping sum of (\ref{eq:lemmaLdiff}), we obtain
\begin{equation}
\begin{aligned}
&\big(\frac{\rho}{4 M_{2}^{2}}-\frac{\gamma+CM_{1}^{2}}{2}-\frac{CM_{1}^2}{\rho} \big) \sum_{\phi = 1}^{\nu} \sum_{k\in\mathcal{A}_{\phi}}\|x_{k}^{\phi+1}-x_{k}^{\phi}\|^2+(\rho+\alpha)\sum_{\phi = 1}^{\nu}\|z^{\phi + 1} - z^{\phi}\|^{2} \\
& ~~~~- \frac{2\rho M_{2}^{4}+1}{2}\sum_{\phi = 1}^{\nu}\sum_{k\in\mathcal{A}_{\phi}}\|z_{k}^{\bar{\phi}_{k}+1}-z_{k}^{\phi}\|^{2}\\
& \leq L(x^1, z^1, \lambda^1)- L(x^{\nu+1}, z^{\nu+1}, \lambda^{\nu+1}) < \infty,
\end{aligned}
\label{eq:Ldifftelescopeinter}
\end{equation}
\normalsize
where the last inequality holds under Assumption \ref{assumption:Lbound}. 

By substituting (\ref{eq:errorbound}) in Lemma (\ref{lemma:errorbound}) into (\ref{eq:Ldifftelescopeinter}), we have
\begin{equation}
\begin{aligned}
&\bigg(\frac{\rho}{4 M_{2}^{2}}-\frac{\gamma+CM_{1}^{2}}{2}-\frac{CM_{1}^2}{\rho} \bigg) \sum_{\phi = 1}^{\nu} \sum_{k\in\mathcal{A}_{\phi}}\|x_{k}^{\phi+1}-x_{k}^{\phi}\|^2 \\
&  +\bigg(\rho+\alpha-\frac{(2\rho M_{2}^{4}+1)(\omega-1)^{2}}{2}\bigg)\sum_{\phi = 1}^{\nu}\|z^{\phi + 1} - z^{\phi}\|^{2} < \infty,
\end{aligned}
\label{eq:Ldifftelescope}
\end{equation}
\normalsize
Then by choosing $\rho$ and $\alpha$ as in (\ref{eq:rhoalpha}), the left-hand-side (LHS) of (\ref{eq:Ldifftelescope}) is positive and increasing with $\nu$. Since the right-hand-side (RHS) of (\ref{eq:Ldifftelescope}) is finite, we must have as $\nu \rightarrow \infty$,
\begin{equation}
x_{k}^{\nu+1} - x_{k}^{\nu} \rightarrow \bm{0}, ~~~~ z_{k}^{\nu+1} - z_{k}^{\nu} \rightarrow \bm{0}, ~~~\forall k.
\label{eq:xzconverge}
\end{equation}
\normalsize
Given (\ref{eq:xzconverge}) and (\ref{eq:lambdabound}), we also have
\begin{equation}
\lambda_{k}^{\nu+1} - \lambda_{k}^{\nu} \rightarrow \bm{0}, ~~~\forall k
\label{eq:lambdaconverge}
\end{equation} 
\normalsize
Since $\mathcal{X}$ and $\mathcal{Z}$ are both compact and $x \in \mathcal{X}$ and $z \in \mathcal{Z}$, and $\lambda$ is bounded by projection, $(\{x_{k}^{\star}\}_{k = 1}^{K}, z^{\star},\{\lambda_{k}^{\star}\}_{k = 1}^{K})$ is bounded and has a limit point. Finally, we show that every limit point of the above sequence is a KKT point of problem (\ref{eq:OPF}); i.e., it satisfies (\ref{eq:KKTall}).

For $k \in \mathcal{A}_{\nu}$, by applying (\ref{eq:lambdaconverge}) to (\ref{eq:lambdaupdateglobal}) and by (\ref{eq:xzconverge}), we obtain
\begin{equation}
A_{k}x_{k}^{\nu+1} - z_{k}^{\nu + 1} \rightarrow \bm{0}, ~~~ k\in \mathcal{A}_{\nu}.
\end{equation}
\normalsize
For $k \notin \mathcal{A}_{\nu}$, let $\bar{\nu}_{k}$ denote the iteration number of region $k$'s last update, then $k \in \mathcal{A}_{\bar{\nu}_{k}}$. Then at iteration $\bar{\nu}_{k}$, we have
\begin{equation}
\lambda_{k}^{\bar{\nu}_{k}+1} = \lambda_{k}^{\bar{\nu}_{k}} + \rho(A_{k}x_{k}^{\bar{\nu}_{k}+1} - z_{k}^{\overline{(\bar{\nu}_{k})}_{k}}),
\end{equation}
\normalsize
where $\overline{(\bar{\nu}_{k})}_{k}$ denotes the iteration of $k$'s update before $\bar{\nu}_{k}$. And since $x_{k}^{\nu+1} = x_{k}^{\nu} = \cdots = x_{k}^{\bar{\nu}_{k}+2} =  x_{k}^{\bar{\nu}_{k}+1} $, and by (\ref{eq:xzconverge}) and (\ref{eq:lambdaconverge}), we have
\begin{equation}
\begin{aligned}
&\|A_{k}x_{k}^{\nu+1}-z_{k}^{\nu+1}\| \\ &= \|A_{k}x_{k}^{\bar{\nu}_{k}+1}-z_{k}^{\nu+1}\| \\
& = \|A_{k}x_{k}^{\bar{\nu}_{k}+1}-z_{k}^{\overline{(\bar{\nu}_{k})}_{k}+1}+z_{k}^{\overline{(\bar{\nu}_{k})}_{k}+1}-z_{k}^{\nu+1}\| \\
&\leq \frac{1}{\rho}\|\lambda_{k}^{\bar{\nu}_{k}+1} - \lambda_{k}^{\bar{\nu}_{k}} \| + \|z_{k}^{\overline{(\bar{\nu}_{k})}_{k}+1}-z_{k}^{\nu+1}\| \rightarrow \bm{0}.
\end{aligned}
\end{equation}
\normalsize
Therefore, we can conclude 
\begin{equation}
A_{k}x_{k}^{\nu+1}-z_{k}^{\nu+1} \rightarrow \bm{0}, \forall k; 
\label{eq:feasibility}
\end{equation}
\normalsize
i.e., the KKT condition (\ref{eq:KKTlambda}) can be satisfied asymptotically.

For any $z_{[ij]}^{\nu+1}$, where $i\in \mathcal{R}_{k}$, $j \in \mathcal{R}_{l}$, $(i,j) \in \mathcal{T}$, and $k \in \mathcal{A}_{\nu}$, the optimality condition of (\ref{eq:zupdateglobal}) yields (\ref{eq:zoptimality}). Since the last three terms on the LHS of (\ref{eq:zoptimality}) will asymptotically converge to $\bm{0}$ due to (\ref{eq:feasibility}) and (\ref{eq:xzconverge}), the KKT condition (\ref{eq:KKTz}) can be satisfied asymptotically. At last, by applying (\ref{eq:xzconverge}) and (\ref{eq:lambdaconverge}) to (\ref{eq:xoptimality}), we obtain KKT condition (\ref{eq:KKTx}). Therefore, we can conclude that $(\{x_{k}^{\nu}\}_{k=1}^{K}, z^{\nu}, \{\lambda_{k}^{\nu}\}_{k=1}^{K})$ are bounded and converge to the set of KKT points of problem (\ref{eq:OPF}).
\end{proof}

\section{Application: AC OPF Problem}
\label{results}
To verify the convergence of the proposed asynchronous ADMM method, we apply Algorithm \ref{algADMMasync} to solve the standard AC OPF problem.
\subsection{Problem Formulation}
The objective of the AC OPF problem is to minimize the total generation cost. The OPF problem is formulated as follows:
\begin{subequations}
\label{centOPF}
\begin{align}
\label{eq1}
\underset{V, P, Q} {\text{minimize}}~~&f({P})=\sum_{i=1}^{n_{b}} \left(a_{i}P_{i}^2+b_{i}P_{i}+c_{i}\right) \\
\label{eqpf}
\text{subject to}~~&P_{i}+jQ_{i}-P_{i}^{\text{load}}-jQ_{i}^{\text{load}}=V_i\sum\limits_{j\in {{\Omega }_{i}}}Y_{ij}^*V_j^*\\
\label{eqPlimit}
&P_{i}^{\min}\leq P_{i}\leq P_{i}^{\max}\\
\label{eqQlimit}
&Q_{i}^{\min}\leq Q_{i}\leq Q_{i}^{\max}\\
\label{eqVlimit}
&V_{i}^{\min}\leq |V_{i}|\leq V_{i}^{\max},
\end{align}
\end{subequations}
\normalsize
for $i=1, \ldots, n_{b}$ where $n_{b}$ is the number of buses. $(a_i,b_i,c_i)$ denote the cost parameters of generator at bus $i$, and $(V_{i}, P_{i}, Q_{i})$ denote the complex voltage, the active and reactive power generation at bus $i$. $Y_{ij}$ is the $ij$-th entry of the line admittance matrix, and $\Omega_i$ is the set of buses connected to bus $i$. This problem is non-convex due to the non-convexity of the AC power flow equations (\ref{eqpf}). We divide the system into regions and use ${x}_{k}$ to denote all the variables in region $k$..
Consequently, constraints (\ref{eqpf}) at the boundary buses are the coupling constraints. We remove such coupling by duplicating the voltages at the boundary buses. For example, assume region $k$ and $l$ are connected via tie line $ij$ with bus $i$ in region $k$ and bus $j$ in region $l$. The voltages at bus $i$ and $j$ are duplicated, and the copies assigned to region $k$ are $V_{i,k}$ and $V_{j,k}$. Similarly, region $l$ is assigned the copies $V_{i,l}$ and $V_{j,l}$. To ensure equivalence with the original problem, constraints $V_{i,k}=V_{i,l}$ and $V_{j,k}=V_{j,l}$ are added to the problem. Then for each tie line $ij$, we introduce the following auxiliary variables \cite{erseghe2015distributed}\cite{guo2016acase} to region $k$: 
\begin{equation}
z_{k,[ij]}^{-}=\beta^{-}(V_{i,k}-V_{j,k}),~~~~~z_{k,[ij]}^{+}=\beta^{+}(V_{i,k}+V_{j,k}),
\label{eq:zV}
\end{equation}
\normalsize
where $\beta^{-} = 2$ and $\beta^{+} = 0.5$ are used in simulations. $\beta^{-}$ should be set to a larger value than $\beta^{+}$ to emphasize on $V_{i}-V_{j}$ which is strongly related to the line flow through line $ij$ \cite{erseghe2015distributed}. 
Similarly, variables $z_{l,[ij]}^{-}=\beta^{-}(V_{i,l}-V_{j,l})$ and $z_{l,[ij]}^{+}=\beta^{+}(V_{i,l}+V_{j,l})$ are introduced to Region $l$. Then $z_{k,[ij]}^{\pm} = z_{l,[ij]}^{\pm}$. Writing all the $z$'s in a compact form, we transform problem (\ref{centOPF}) to the desired formulation (\ref{eq:OPF}). As the feasible regions of the OPF problem are smooth compact manifolds \cite{chiang2017feasible}, assumptions \ref{assumption:boundeddelay} to \ref{assumption:Lbound} can be satisfied.

\label{simulation}
\subsection{Experiment Setup}
The simulations are conducted using two IEEE standard test systems and two large-scale transmission networks. The system configuration and the parameter settings are given in Table \ref{table:asyncsystemconfig}. The partitions of the systems are derived using the partitioning approach proposed in \cite{guointelligent} that reduces the coupling among regions \cite{guo2016acase}. A ``flat" start initializes $x$ to be at the median of its upper and lower bounds, while a ``warm" start is a feasible solution to the power flow equations (\ref{eqpf}).

Algorithm \ref{algADMMasync} is conducted in Matlab with $p$ set to 0.1 which simulates the worst case where each worker is allowed to perform a local update with one arrived neighbor. The stopping criterion is that the maximum residue ($\max\{\Gamma_k\}, \forall k$) and constraint mismatch are both smaller than $10^{-3}$ p.u. We use number of average local iterations and the execution time to measure the performance of Algorithm \ref{algADMMasync}. The execution time records the total time Algorithm \ref{algADMMasync} takes until convergence including the computation time (measured by CPU time) and the waiting time for neighbors. Here, the waiting time also includes the communication delay, which is estimated by assuming that fiber optical communications is used. Therefore, passing message from one worker to the other usually takes a couple of milliseconds, which is very small compared to local computation time. 
\begin{table}[tb]
\caption{Test system configuration and parameter setting for asynchronous ADMM.}
\vspace{0cm}
\centering
\begin{tabular}{ p{2.2cm}<{\centering}|  p{2.2cm}<{\centering} | p{2.2cm}<{\centering}| p{2.2cm}<{\centering}|p{3.5cm}<{\centering}}
\toprule
system&IEEE-30&IEEE-118&Polish&Great Britain (GB)\\
\midrule
buses&30&118&2383&2224\\
regions&4&8&40&70\\
initialization&flat&flat&warm&warm\\
$\rho$&$5\times 10^4$&$5\times 10^5$&$10^8$&$10^9$\\
\bottomrule
\end{tabular}
\vspace{-0.2cm}
\label{table:asyncsystemconfig}
\end{table}

\begin{figure*}[t]
\setlength{\abovecaptionskip}{0.1cm} 
\centering
\captionsetup[subfigure]{captionskip=-0.00cm}
\subfloat[IEEE-30 Iteration]
{
\label{fig:30iterasync}
\includegraphics[trim = 0mm 0mm 10mm 0mm, clip=true,width=3.8cm]{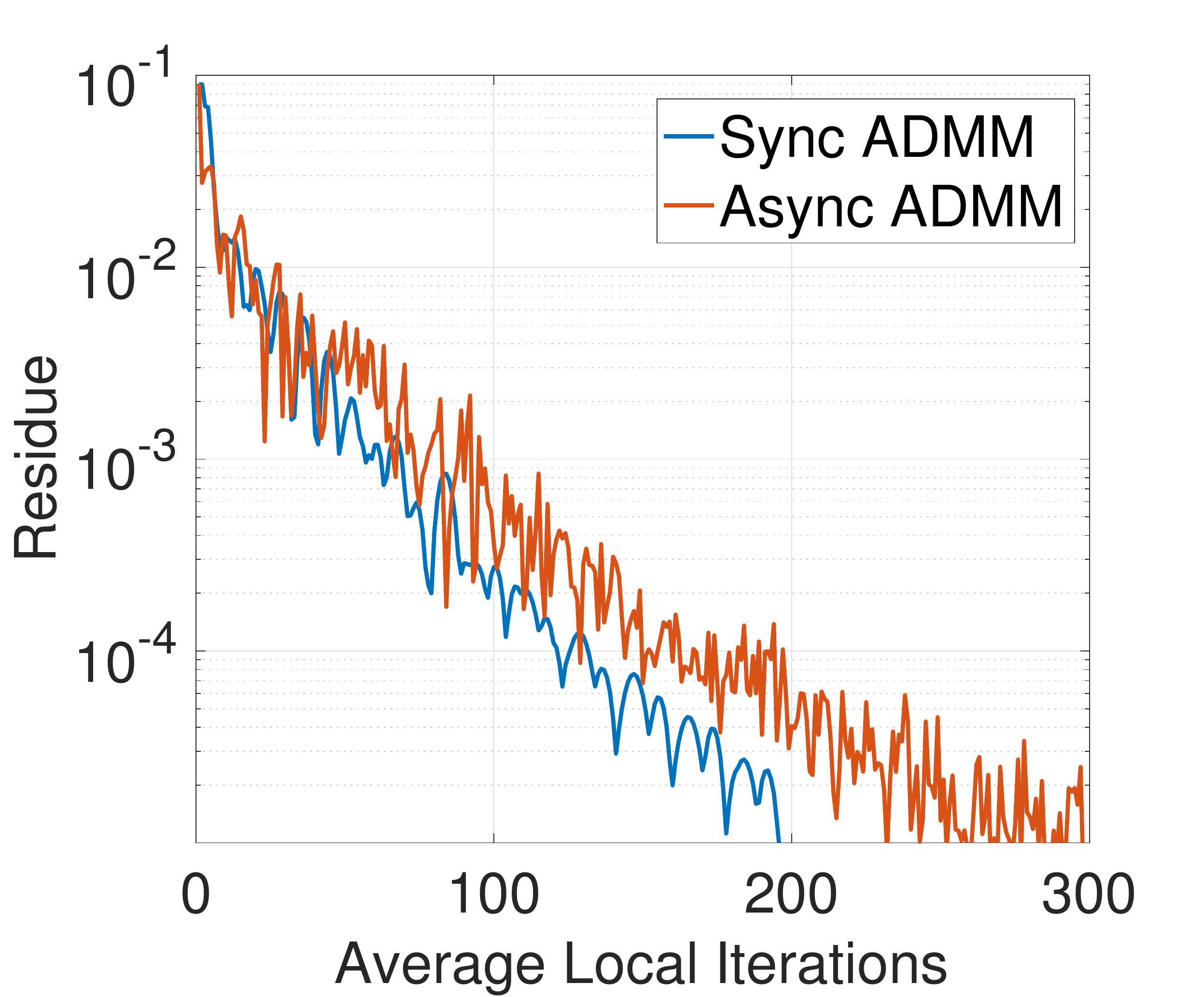} 
}
\hspace{0cm}
\subfloat[IEEE-118 Iteration]{
\label{fig:118iterasync}
\includegraphics[trim = 0mm 0mm 10mm 0mm, clip=true,width=3.8cm]{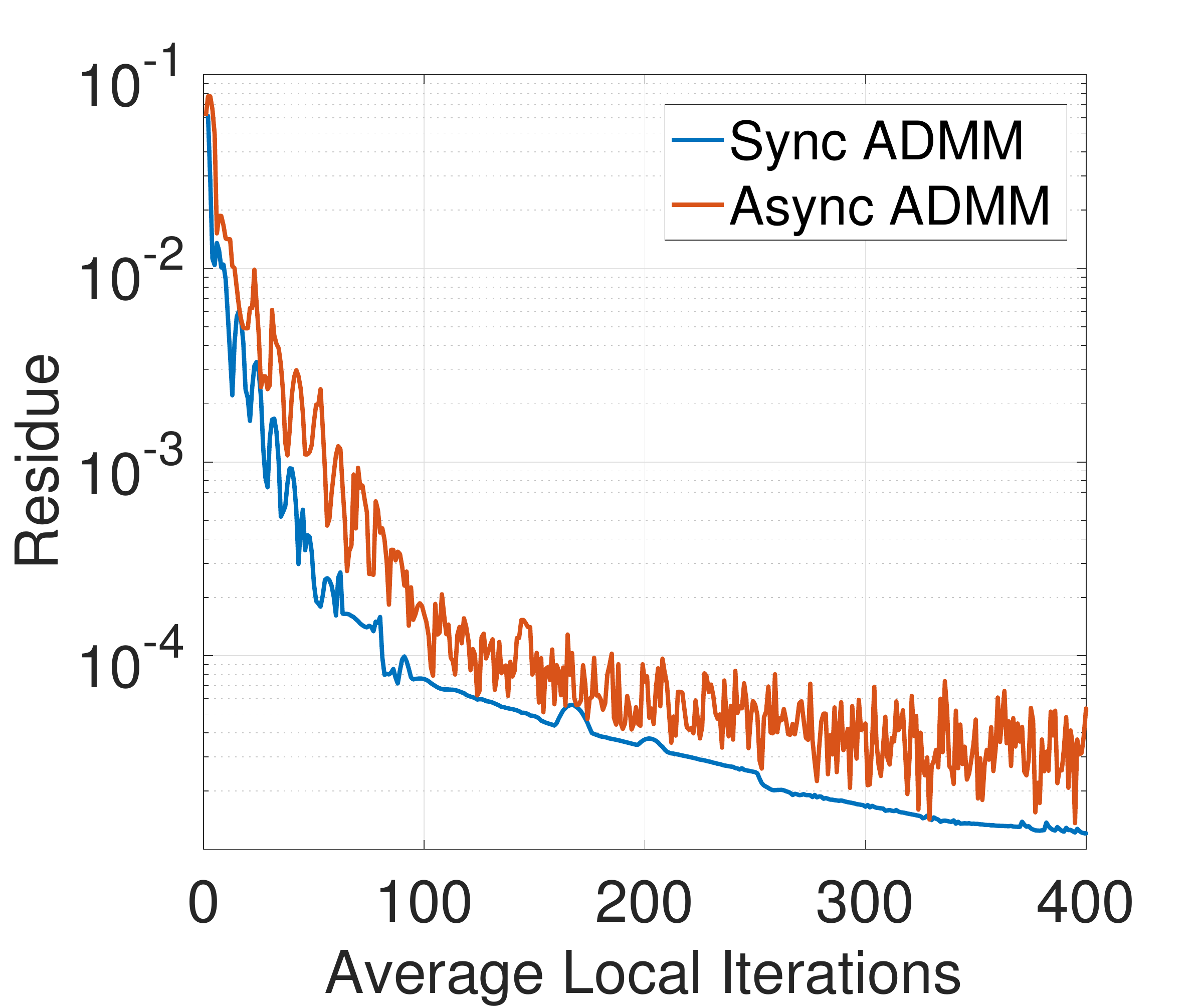} 
}
\hspace{0cm}
\subfloat[Polish Iteration]{
\label{fig:polishiterasync}
\includegraphics[trim = 0mm 0mm 10mm 0mm, clip=true,width=3.8cm]{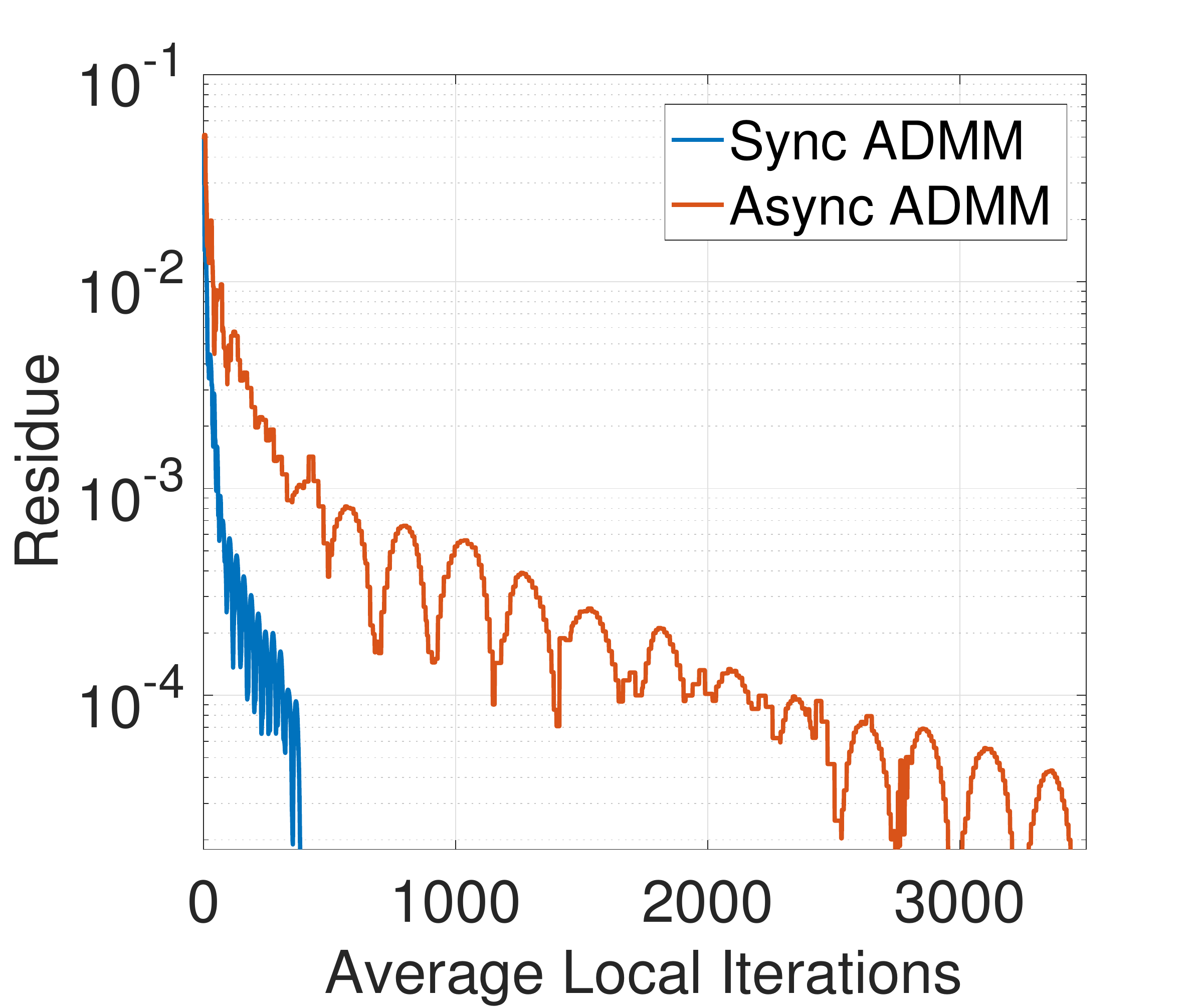} 
}
\hspace{0cm}
\subfloat[GB Iteration]{
\label{fig:gbiterasync}
\includegraphics[trim = 0mm 0mm 10mm 0mm, clip=true,width=3.8cm]{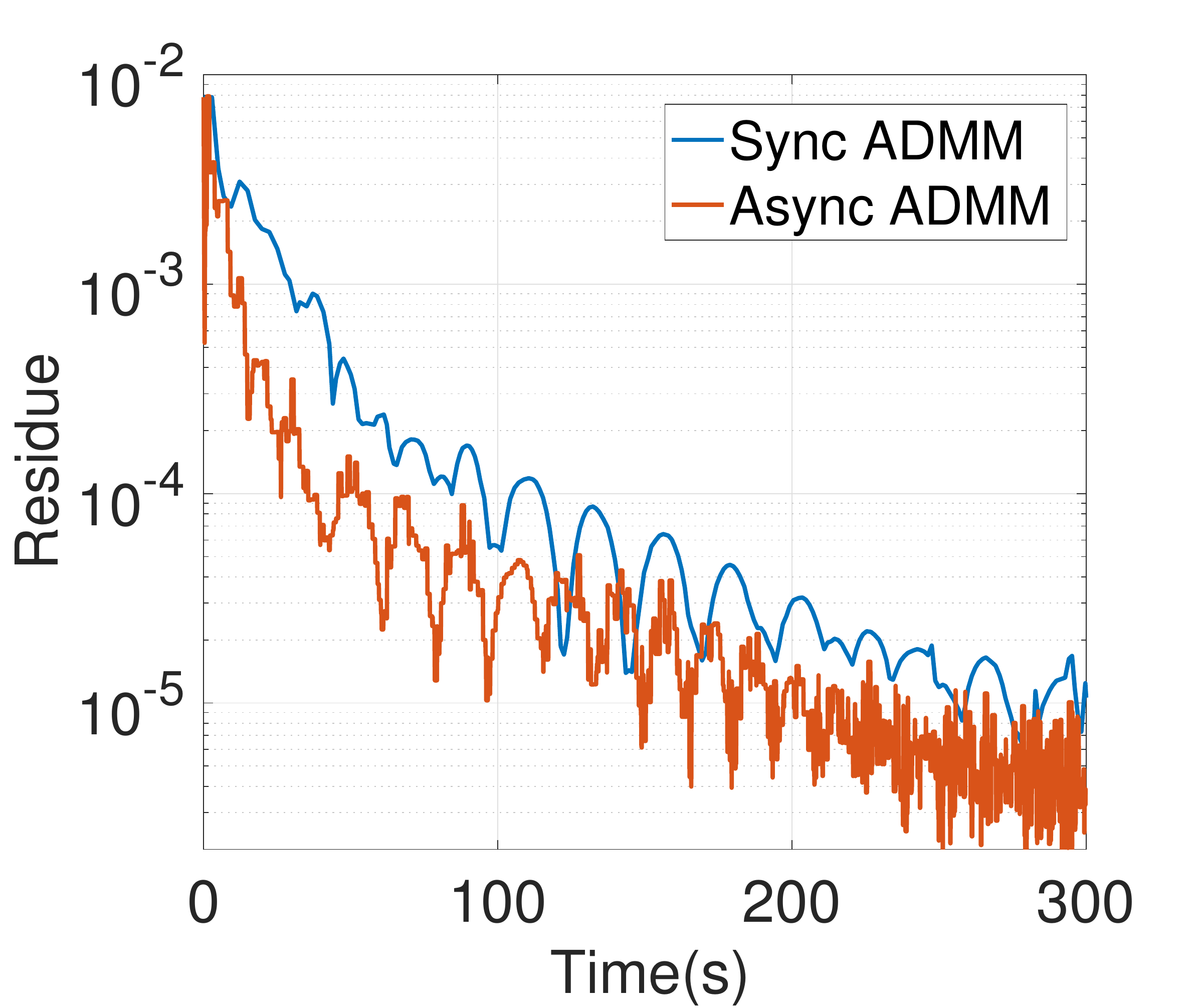} 
}\\

\subfloat[IEEE-30 Time]{
\label{fig:30timeasync}
\includegraphics[trim = 0mm 0mm 10mm 0mm,clip=true, width=3.8cm]{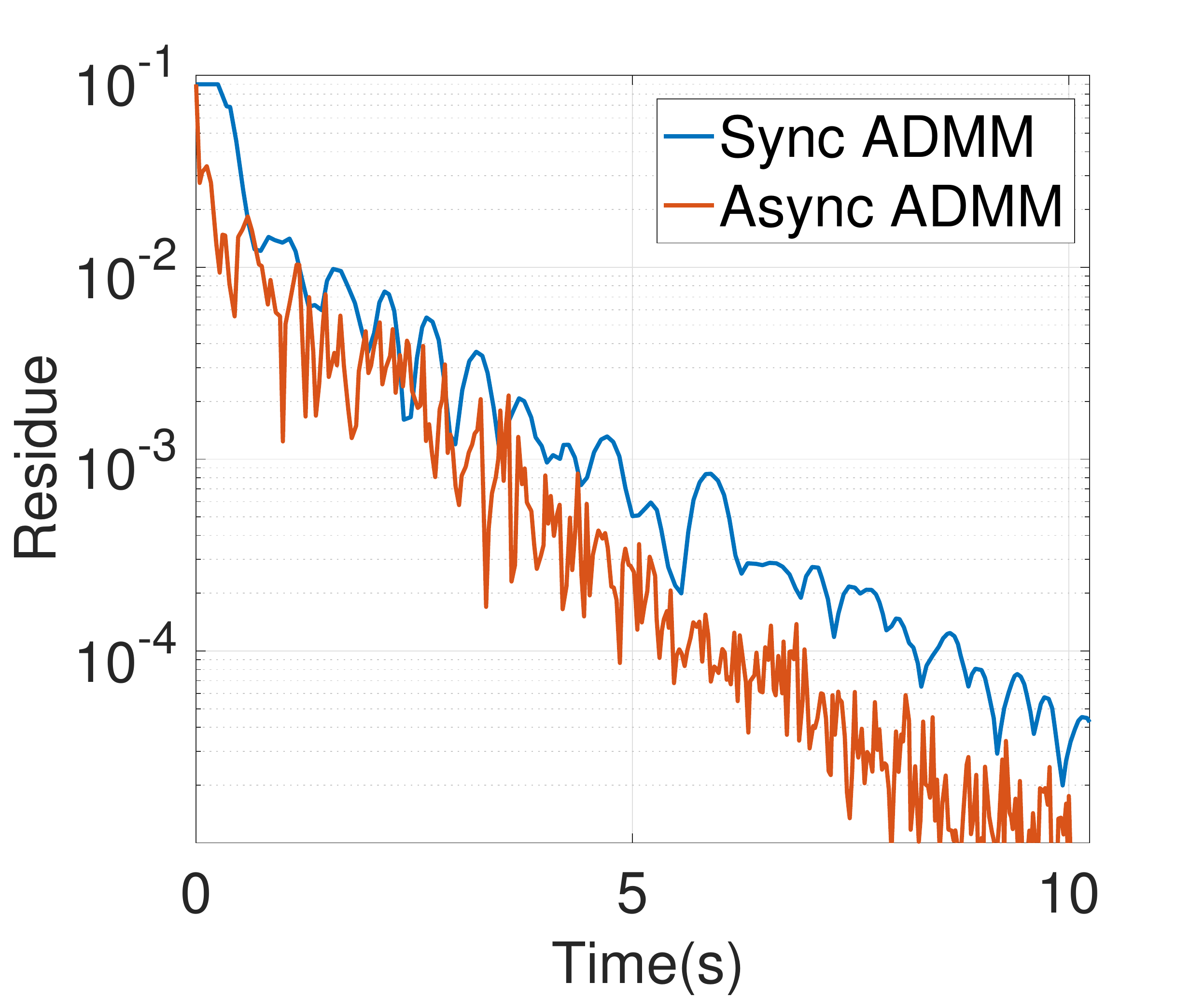} 
}
\hspace{0cm}
\subfloat[IEEE-118 Time]
{
\label{fig:118timeasync}
\includegraphics[trim = 0mm 0mm 10mm 0mm, clip=true,width=3.8cm]{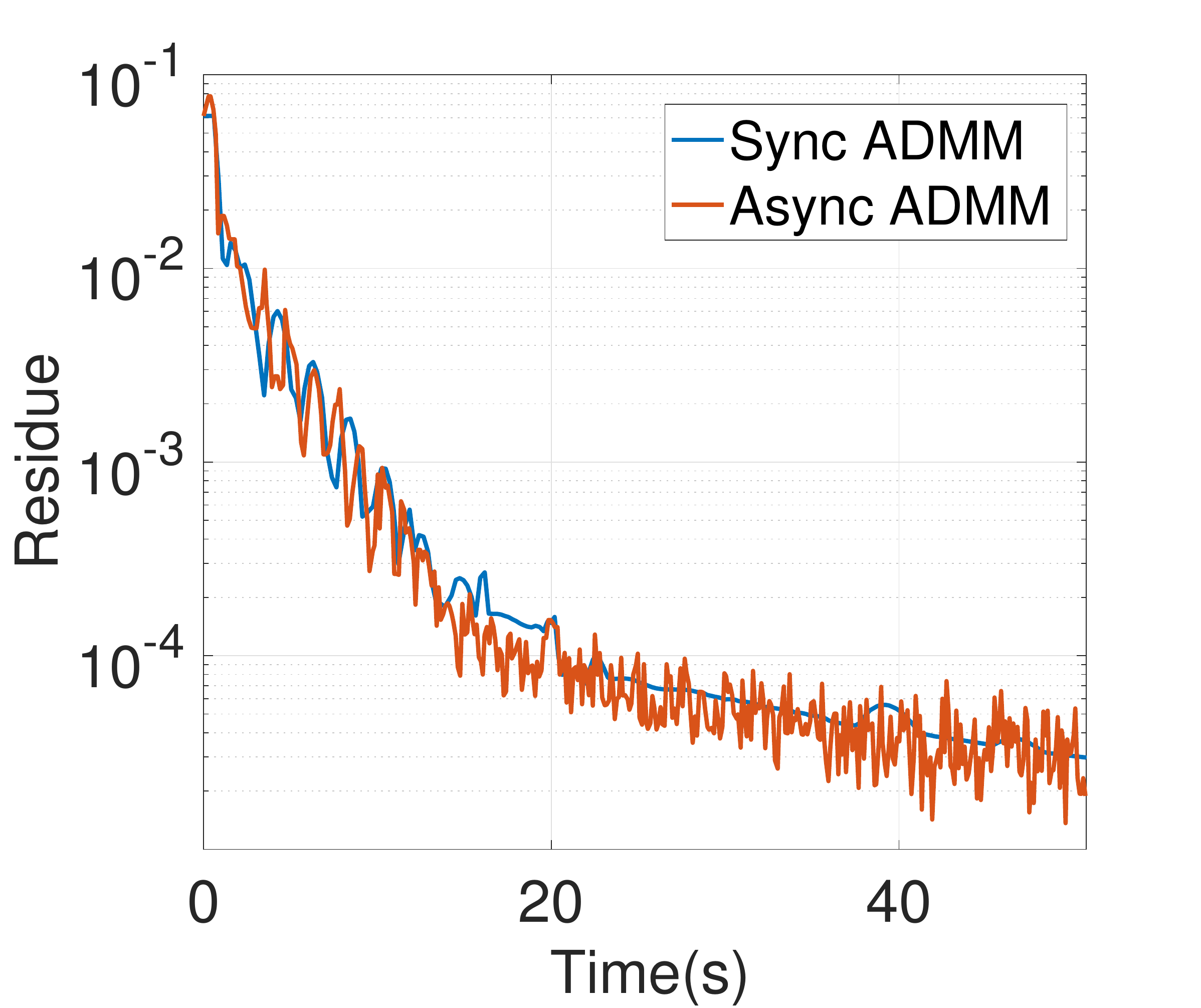} 
}
\hspace{0cm}
\subfloat[Polish Time]
{
\label{fig:polishtimeasync}
\includegraphics[trim = 0mm 0mm 10mm 0mm, clip=true,width=3.8cm]{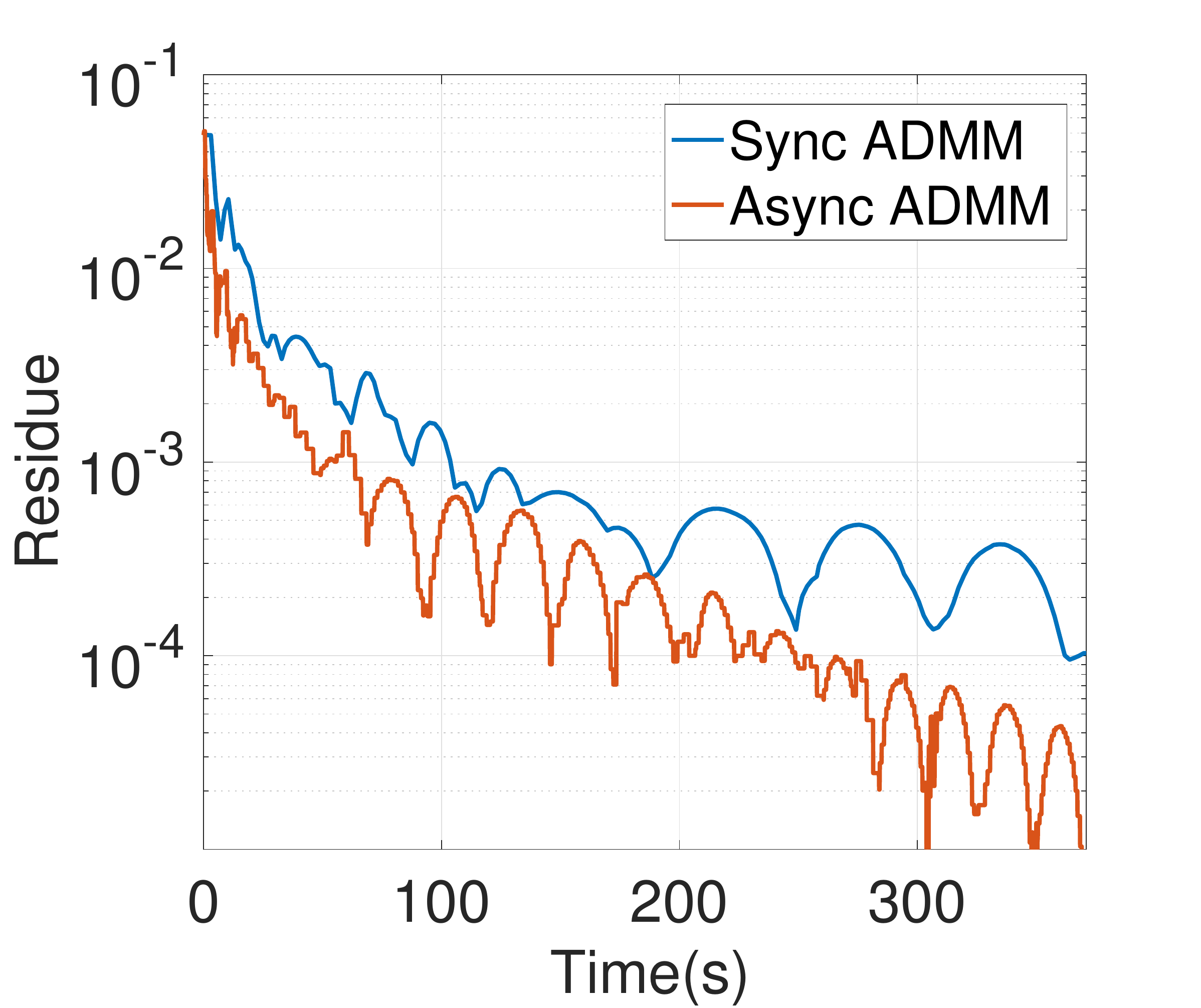} 
}
\hspace{0cm}
\subfloat[GB Time]
{
\label{fig:gbtimeasync}
\includegraphics[trim = 0mm 0mm 10mm 0mm, clip=true,width=3.8cm]{GBtimeasync-eps-converted-to.pdf} 
}
\caption{Convergence of residue of synchronous and asynchronous ADMM.}
\vspace{-0.2cm}
\label{fig:asyncADMMconvsystems}
\end{figure*}
\subsection{Numerical Results}
\begin{figure}[bt]
\setlength{\abovecaptionskip}{0cm} 
\centering
\includegraphics[trim = 10mm 0mm 0mm 0mm, clip=true,width=13cm]{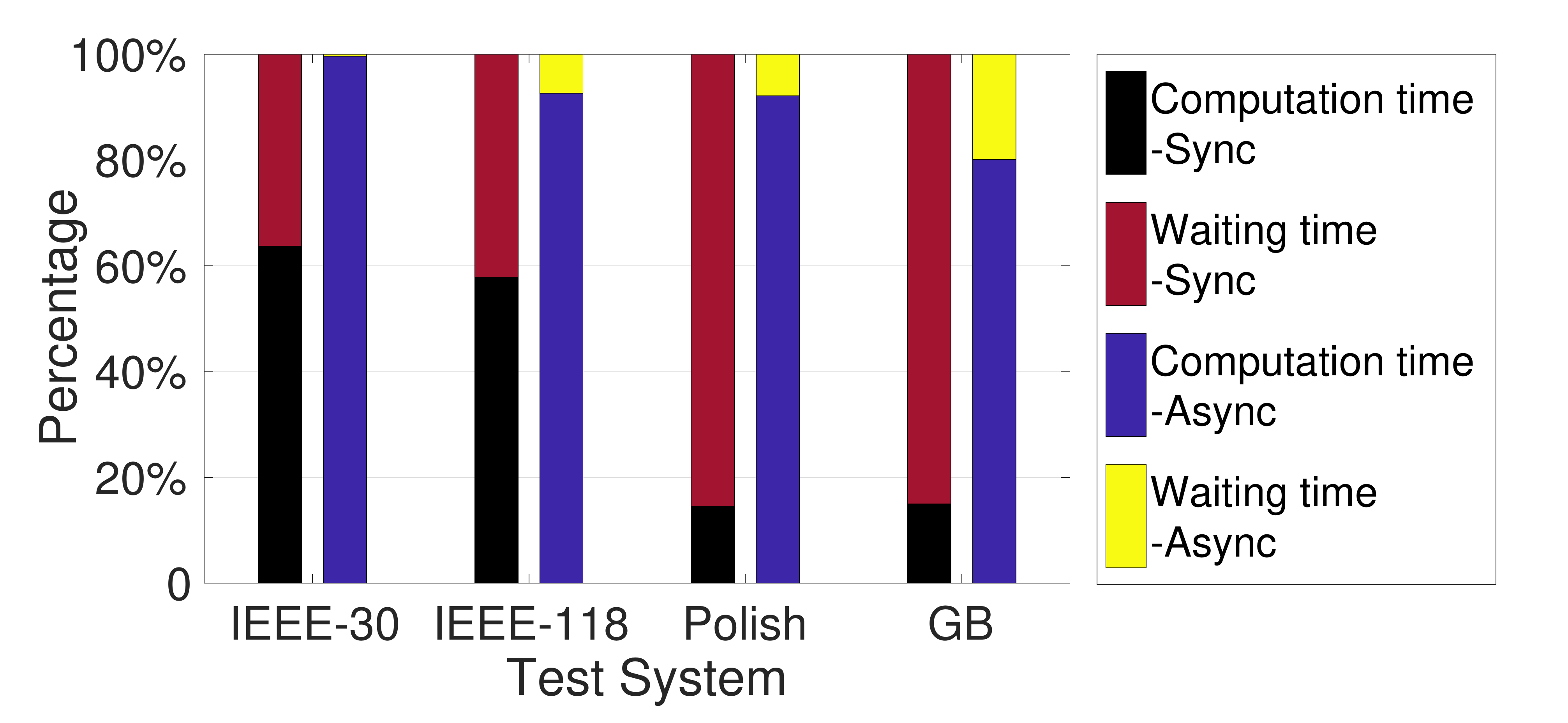} 
\caption{Computation/waiting time for synchronous ADMM and asynchronous ADMM on different test systems.}
\vspace{-0.5cm}
\label{fig:twaitcompasync}
\end{figure}

Figure \ref{fig:asyncADMMconvsystems} shows the convergence of the maximum residue of the proposed asynchronous ADMM and the standard synchronous ADMM on solving the OPF problem for the considered four test systems. We set $\alpha$ to zero for this experiment and will show later that using a large $\alpha$ is not necessary for the convergence of asynchronous ADMM. As expected, synchronous ADMM takes fewer iterations to converge, especially on large-scale systems. However, due to the waiting time for the slowest worker at each iteration, synchronous ADMM can be slower than asynchronous ADMM especially on large-scale networks. Figure \ref{fig:twaitcompasync} illustrates the percentage of the average computation and waiting time experienced by all workers. It is clearly shown that a lot of time is wasted on waiting for all neighbors using a synchronous scheme. 

To measure the optimality of the solution found by asynchronous ADMM, we also calculate the gap in the objective value achieved by synchronous and asynchronous ADMM with respect to the objective value obtained by a centralized method when the maximum residue $\Gamma_k$ of all workers reaches $10^{-3}$. This gap is shown in Table \ref{table:gapasync} which are fairly small for all the considered systems with both schemes. A solution can be considered of good quality if this gap is smaller than $1\%$. Surprisingly, asynchronous ADMM achieves a slightly smaller gap compared with synchronous ADMM. This is due to the fact that in asynchronous ADMM a worker uses the most updated information of its neighbors in a more timely manner than in the synchronous case and updates its local variables more frequently, which, as a trade-off, results in more local iterations especially on large systems. Note that for both schemes, this gap should asymptotically approach zero. However, for many engineering applications, only a mild level of accuracy is needed. Therefore ADMM is usually terminated when it reaches the required accuracy with a suboptimal solution. These results indeed validate that asynchronous ADMM could find a good solution faster than its synchronous counterpart and is more fault-tolerable for delayed or missing information. But we should mention the large number of iterations taken by asynchronous ADMM incurs more communications load because updated information is sent at each iteration. And this might raise more requirements on the communications system used for deploying distributed algorithms, which will subject to our future studies.

\begin{table}[tb]
\caption{Gap of objective function achieved by synchronous and asynchronous ADMM.}
\vspace{0cm}
\centering
\begin{tabular}{ p{2.2cm}<{\centering}|  p{2.2cm}<{\centering} | p{2.2cm}<{\centering}| p{2.2cm}<{\centering}|p{3.5cm}<{\centering}}
\toprule
scheme&IEEE-30&IEEE-118&Polish&Great Britain (GB)\\
\midrule
synchronous&0.025\%&0.122\%&0.060\%&0.575\%\\
asynchronous &0.005\%&0.098\%&0.031\%&0.416\%\\
\bottomrule
\end{tabular}
\label{table:gapasync}
\end{table}

\begin{figure}[t]
\setlength{\abovecaptionskip}{0.0cm} 
\centering
\captionsetup[subfigure]{captionskip=-0.00cm}
\subfloat[Impact of $\alpha$]
{
\label{fig:30asyncalpha}
\includegraphics[trim = 0mm 0mm 10mm 0mm, clip=true,width=7.2cm]{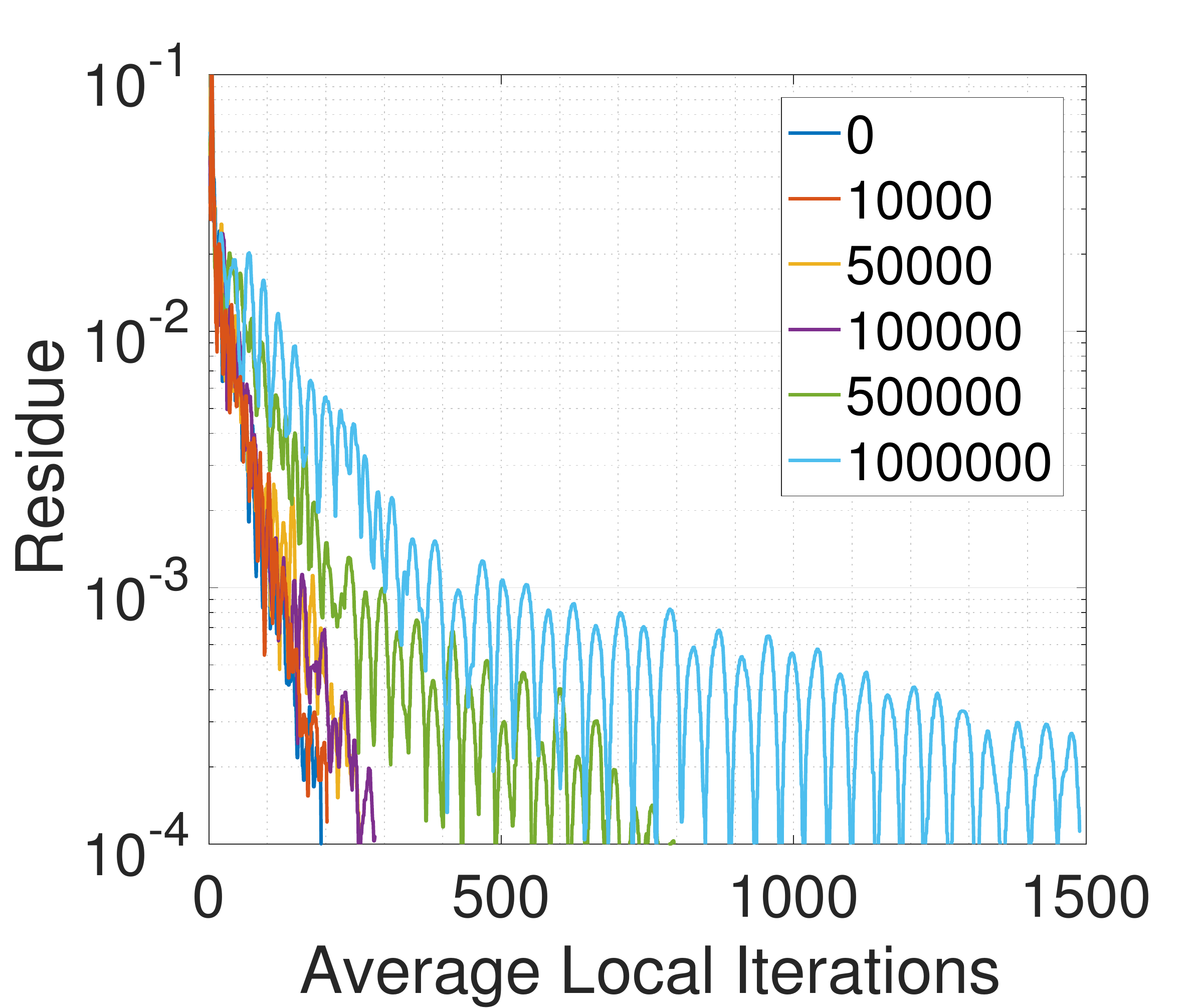} 
}
\hspace{-0.2cm}
\subfloat[Impact of $\rho$]{
\label{fig:30asyncrho}
\includegraphics[trim = 0mm 0mm 10mm 0mm, clip=true,width=7.2cm]{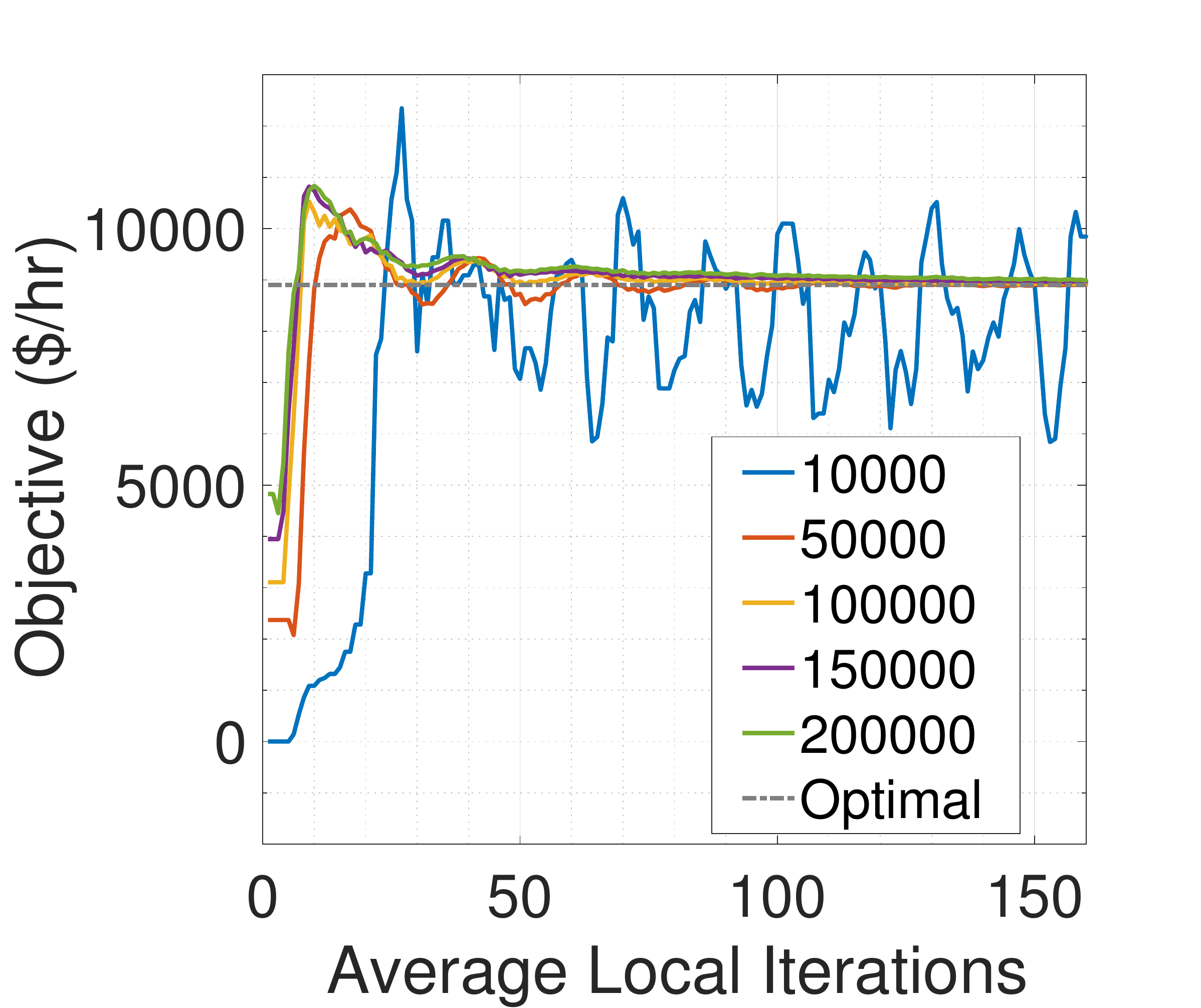} 
}
\caption{Impact of parameters on the convergence of asynchronous ADMM.}
\vspace{-0.2cm}
\label{fig:asyncADMMpar}
\end{figure}
Now, we evaluate the impact of parameter values on the performance of asynchronous ADMM. Interestingly, as shown in Fig. \ref{fig:30asyncalpha}, even though a large $\alpha$ is a sufficient condition for asynchronous ADMM to converge, it is not necessary in practice. This observation is consistent with the observation made in \cite{7423789} as the proof there and our proof are both derived for the worst case. In fact, the purpose of using $\alpha$ is to make sure that the Augmented Lagrangian (\ref{eq:Aug}) decreases at each iteration which may not be the case due to the existence of the term $\rho\sum_{k\in\mathcal{A}_{\nu}}({z}_{k}^{\bar{\nu}_{k}+1} - {z}_{k}^{\nu})^{\top}(A_{k}x_{k}^{\nu+1}-A_{k}x_{k}^{\nu})$ in (\ref{eq:Ldiffinter}). However, as $z_k$ is generally a value between $A_k x_k$ and $A_l x_l, \forall l \in \mathcal{N}_k$ that worker $k$ and $l$ try to approach, $({z}_{k}^{\bar{\nu}_{k}+1} - {z}_{k}^{\nu})^{\top}(A_{k}x_{k}^{\nu+1}-A_{k}x_{k}^{\nu})$ is likely to be a negative value, which makes $\alpha$ unnecessary. In fact, as shown in Fig. \ref{fig:30asyncalpha}, a large $\alpha$ will slow down the convergence as the proximal term in (\ref{eq:zupdatelocal}) forces local updates to take very small steps. Finally, Fig. \ref{fig:30asyncrho} shows that a large $\rho$ is indeed necessary for the convergence of asynchronous ADMM. With a larger $\rho$, asynchronous ADMM tends to stabilize around the final solution more quickly, which, however, may lead to a slightly less optimal solution.

\section{Conclusions and Future Works}
\label{conclusion}
This paper proposes an asynchronous distributed ADMM approach to solve the non-convex optimization problem given in (\ref{eq:intro}). The proposed method does not require any centralized coordination and allows any worker to perform local updates with information received from a subset of its physically connected neighbors. We provided the sufficient conditions under which asynchronous ADMM asymptotically satisfies the first-order optimality conditions. Through the application of the proposed approach to solve the AC OPF problem for several power systems, we demonstrated that asynchronous ADMM could be more efficient than its synchronous counterpart and therefore is more suitable and scalable for distributed optimization applications in large-scale systems.

While the results acquired in this paper provides the theoretical foundation for studies on asynchronous distributed optimization, there are many practical issues that need to be addressed. For example, the presented algorithmic framework does not include the model of communication delay, which may have a strong impact on the convergence of asynchronous distributed methods. Moreover, it is also important to define good system partitions and choose proper algorithm parameters. Therefore, in the future, we plan to investigate how those different factors affect the convergence speed of the proposed asynchronous distributed optimization scheme.
\section*{Acknowledgment}
The authors would like to thank ABB for the financial support and Dr. Xiaoming Feng for his invaluable inputs.
    \vspace{-0.3cm}
\appendix
    \vspace{-0.4cm}
\section{Proof of Lemma~\ref{lemma:Ldifference}}
\vspace{-0.1cm}
\begin{proof}[Proof of Lemma~\ref{lemma:Ldifference}]
\begin{equation}
\begin{aligned}
&L(x^{\nu+1}, z^{\nu+1}, \lambda^{\nu+1})-L(x^{\nu}, z^{\nu}, \lambda^{\nu})\\&=L(x^{\nu+1}, z^{\nu}, \lambda^{\nu})-L(x^{\nu}, z^{\nu}, \lambda^{\nu})\\
&~~~+L(x^{\nu+1}, z^{\nu}, \lambda^{\nu+1})-L(x^{\nu+1}, z^{\nu}, \lambda^{\nu})\\&~~~+L(x^{\nu+1}, z^{\nu+1}, \lambda^{\nu+1})-L(x^{\nu+1}, z^{\nu}, \lambda^{\nu+1}),
\end{aligned}
\label{eq:Ldifference}
\normalsize
\end{equation}
We bound the three pairs of differences on the RHS of (\ref{eq:Ldifference}) as follows. 
First, due to the optimality of $x_{k}^{\nu+1}$ in (\ref{eq:xupdateglobal}), we introduce the general subgradient $d_{k}^{\nu + 1}$:
\begin{equation}
d_{k}^{\nu+1}:= -\big(A_{k}^{\top}\lambda_{k}^\nu+\rho A_{k}^{\top}(A_{k}{x}_{k}^{\nu+1}-{z}_{k}^{\bar{\nu}_{k}+1})\big)\in \partial F_{k}(x_{k}^{\nu+1})
\label{eq:d}
\end{equation}
\normalsize
Then, since only $x_{k}$ for $k\in\mathcal{A}_{\nu}$ is updated, we have
\begin{equation}
\begin{aligned}
&L(x^{\nu}, z^{\nu}, \lambda^{\nu}) - L(x^{\nu+1}, z^{\nu}, \lambda^{\nu}) \\
& = \sum_{k\in\mathcal{A}_{\nu}} \Big(F_{k}({x}_{k}^{\nu}) - F_{k}({x}_{k}^{\nu+1})+{\lambda}_{k}^{\nu\top}A_{k}({x}_{k}^{\nu}-x_{k}^{\nu+1})+\frac{{\rho}}{2}\|A_{k}{x}_{k}^{\nu}-{z}_{k}^{\nu}\|^{2}-\frac{{\rho}}{2}\|A_{k}{x}_{k}^{\nu+1}-{z}_{k}^{\nu}\|^{2} \Big)\\
&= \sum_{k\in\mathcal{A}_{\nu}} \Big(F_{k}({x}_{k}^{\nu}) - F_{k}({x}_{k}^{\nu+1})+{\lambda}_{k}^{\nu\top}A_{k}({x}_{k}^{\nu}-x_{k}^{\nu+1}) + \frac{\rho}{2}\|A_{k}x_{k}^{\nu}-A_{k}x_{k}^{\nu+1}\|^2 \\
&~~~~+\rho \langle A_{k}{x}_{k}^{\nu+1}-{z}_{k}^{\bar{\nu}_{k}+1} + {z}_{k}^{\bar{\nu}_{k}+1} - {z}_{k}^{\nu}, A_{k}x_{k}^{\nu}-A_{k}x_{k}^{\nu+1} \rangle\Big)\\
&=\sum_{k\in\mathcal{A}_{\nu}} \Big(F_{k}({x}_{k}^{\nu}) - F_{k}({x}_{k}^{\nu+1})+ \frac{\rho}{2}\|A_{k}x_{k}^{\nu}-A_{k}x_{k}^{\nu+1}\|^2 + \langle A_{k}^{\top}\lambda_{k}^\nu + \rho A_{k}^{\top}(A_{k}{x}_{k}^{\nu+1}-{z}_{k}^{\bar{\nu}_{k}+1}), x_{k}^{\nu}-x_{k}^{\nu+1} \rangle\\&~~~~+\rho({z}_{k}^{\bar{\nu}_{k}+1} - {z}_{k}^{\nu})^{\top}(A_{k}x_{k}^{\nu}-A_{k}x_{k}^{\nu+1})\Big)\\
&= \sum_{k\in\mathcal{A}_{\nu}} \Big(F_{k}({x}_{k}^{\nu}) - F_{k}({x}_{k}^{\nu+1}) - \langle d_{k}^{\nu+1}, x_{k}^{\nu}-x_{k}^{\nu+1} \rangle +\rho({z}_{k}^{\bar{\nu}_{k}+1} - {z}_{k}^{\nu})^{\top}(A_{k}x_{k}^{\nu}-A_{k}x_{k}^{\nu+1})\Big)\\
&~~~~+ \frac{\rho}{2}\|A_{k}x_{k}^{\nu}-A_{k}x_{k}^{\nu+1}\|^2\\
& \geq - \frac{\gamma}{2}\sum_{k\in\mathcal{A}_{\nu}}\|x_{k}^{\nu+1}-x_{k}^{\nu}\|^2 + \frac{\rho}{2}\sum_{k\in\mathcal{A}_{\nu}} \|A_{k}x_{k}^{\nu}-A_{k}x_{k}^{\nu+1}\|^2 + \rho\sum_{k\in\mathcal{A}_{\nu}} ({z}_{k}^{\bar{\nu}_{k}+1} - {z}_{k}^{\nu})^{\top}(A_{k}x_{k}^{\nu}-A_{k}x_{k}^{\nu+1}),
\end{aligned}
\label{eq:Ldiffx}
\end{equation}
\normalsize
where the second equality follows from the cosine rule: \small\mbox{$\|b+c\|^{2}-\|a+c\|^{2}=\|b-a\|^{2}+2\langle a+c, b-a \rangle$}\normalsize with \small$a=A_{k}x_{k}^{\nu+1}$, $b=A_{k}x_{k}^{\nu}$\normalsize, and \small$c=-z_{k}^{\nu}$\normalsize, and the fourth equality is due to the definition of \small$d_{k}^{\nu+1}$\normalsize in (\ref{eq:d}). The last inequality is derived from Definition \ref{def:proxregular} under Assumption \ref{assumption:xfeasible}. Then, by Assumption \ref{assumption:Axmapping}, we have
\begin{equation}
\begin{aligned}
& L(x^{\nu+1}, z^{\nu}, \lambda^{\nu}) -L(x^{\nu}, z^{\nu}, \lambda^{\nu}) \\
& \leq \frac{\gamma}{2}\sum_{k\in\mathcal{A}_{\nu}}\|x_{k}^{\nu+1}-x_{k}^{\nu}\|^2 - \frac{\rho}{2}\sum_{k\in\mathcal{A}_{\nu}}\|A_{k}x_{k}^{\nu}-A_{k}x_{k}^{\nu+1}\|^2-\rho\sum_{k\in\mathcal{A}_{\nu}}({z}_{k}^{\bar{\nu}_{k}+1} - {z}_{k}^{\nu})^{\top}(A_{k}x_{k}^{\nu}-A_{k}x_{k}^{\nu+1})\\
& \leq \frac{\gamma M_{2}^{2} - \rho}{2 M_{2}^{2}} \sum_{k\in\mathcal{A}_{\nu}}\|x_{k}^{\nu+1}-x_{k}^{\nu}\|^2+\rho\sum_{k\in\mathcal{A}_{\nu}}({z}_{k}^{\bar{\nu}_{k}+1} - {z}_{k}^{\nu})^{\top}(A_{k}x_{k}^{\nu+1}-A_{k}x_{k}^{\nu}).
\end{aligned}
\label{eq:Ldiffxfinal}
\end{equation}
\normalsize
Secondly, for the $\lambda$-update, we have\begin{equation}
\begin{aligned}
&L(x^{\nu+1}, z^{\nu}, \lambda^{\nu+1})-L(x^{\nu+1}, z^{\nu}, \lambda^{\nu})\\
&= \sum_{k\in\mathcal{A}_{\nu}}(\lambda_{k}^{\nu+1} - \lambda_{k}^{\nu})^{\top}(A_{k}x_{k}^{\nu+1} - z_{k}^{\nu})\\
&=\sum_{k\in\mathcal{A}_{\nu}}(\lambda_{k}^{\nu+1} - \lambda_{k}^{\nu})^{\top}(A_{k}x_{k}^{\nu+1} - z_{k}^{\bar{\nu}_{k}+1})+\sum_{k\in\mathcal{A}_{\nu}}(\lambda_{k}^{\nu+1} - \lambda_{k}^{\nu})^{\top}(z_{k}^{\bar{\nu}_{k}+1}- z_{k}^{\nu})\\
&=\frac{1}{\rho}\sum_{k\in\mathcal{A}_{\nu}}\|\lambda_{k}^{\nu+1} - \lambda_{k}^{\nu}\|^{2} + \sum_{k\in\mathcal{A}_{\nu}}(\lambda_{k}^{\nu+1} - \lambda_{k}^{\nu})^{\top}(z_{k}^{\bar{\nu}_{k}+1}- z_{k}^{\nu}).
\end{aligned}
\label{eq:Ldifflambda}
\end{equation}
\normalsize
The first equality is due to the fact that \small$\lambda_{k}^{\nu+1} = \lambda_{k}^{\nu}, \forall k \notin \mathcal{A}_{\nu}$\normalsize, and the last equality is obtained by applying (\ref{eq:lambdaupdateglobal}) for \small$k \in \mathcal{A}_{\nu}$\normalsize.

Thirdly, for the $z$-update, we first look at the difference caused by updating $z_{k,l}$ (or $z_{l,k}$ for one pair of neighboring workers $k$ and $l$. By Definition \ref{def:strongconvexity}.1, (\ref{eq:Lsingz}) is strongly convex with respect to (w.r.t.) $z_{k,l}$ with modulus \small$2\rho + \alpha$\normalsize. Then by Definition \ref{def:strongconvexity}.2, we have
\begin{equation}
\begin{aligned}
&L'(x^{\nu+1}, z_{k,l}^{\nu}, \lambda^{\nu + 1}) - L'(x^{\nu+1}, z_{k,l}^{\nu+1}, \lambda^{\nu + 1})\\
&=\bigg(-(\lambda_{k,l}^{\nu + 1} +  \lambda_{l,k}^{\nu + 1})z_{k,l}^{\nu} + \frac{\rho}{2}\|A_{k,l}x_{k}^{\nu+1}-z_{k,l}^{\nu}\|^{2}+ \frac{\rho}{2}\|A_{l,k}x_{l}^{\nu+1}-z_{k,l}^{\nu}\|^{2}+\frac{\alpha}{2}\|z_{k,l}^{\nu}-z_{k,l}^{\nu}\|^{2} \bigg)\\&~~~~-\bigg(-(\lambda_{k,l}^{\nu + 1} +  \lambda_{l,k}^{\nu + 1})z_{k,l}^{\nu + 1} + \frac{\rho}{2}\|A_{k,l}x_{k}^{\nu+1}-z_{k,l}^{\nu + 1}\|^{2}+ \frac{\rho}{2}\|A_{l,k}x_{l}^{\nu+1}-z_{k,l}^{\nu + 1}\|^{2} +\frac{\alpha}{2}\|z_{k,l}^{\nu+1}-z_{k,l}^{\nu}\|^{2}\bigg)\\
&\geq \bigg( \frac{\partial L'(x^{\nu + 1}, z_{k,l}^{\nu + 1}, \lambda^{\nu + 1})}{ \partial z_{k,l}^{\nu + 1}} (z_{k,l}^{\nu} - z_{k,l}^{\nu + 1}) + (\rho+\frac{\alpha}{2})\|z_{k,l}^{\nu } - z_{k,l}^{\nu + 1}\|^{2}  \bigg) \\
&= (\rho+\frac{\alpha}{2})\|z_{k,l}^{\nu } - z_{k,l}^{\nu + 1}\|^{2} 
\end{aligned}
\label{eq:Ldiffsingz}
\end{equation}
\normalsize
The last equality in (\ref{eq:Ldiffsingz}) holds because \small$\partial L'(z_{k,l}^{\nu + 1}) = 0~$\normalsize due to the optimality condition of (\ref{eq:zupdateglobal}). For any \small$z_{k,l}\in z$\normalsize, we have 
\begin{equation}
\begin{aligned}
&L(x^{\nu+1}, z_{k,l}^{\nu + 1}, \lambda^{\nu + 1}) - L(x^{\nu+1}, z_{k,l}^{\nu}, \lambda^{\nu + 1}) \leq -(\rho+\alpha)\|z_{k,l}^{\nu + 1} - z_{k,l}^{\nu}\|^{2}.
\end{aligned}
\label{eq:Ldiffanyz}
\end{equation}
\normalsize
Summing up the LHS of (\ref{eq:Ldiffanyz}), we have
\begin{equation}
\begin{aligned}
&L(x^{\nu+1}, z^{\nu + 1}, \lambda^{\nu + 1}) - L(x^{\nu+1}, z^{\nu}, \lambda^{\nu + 1})\leq -(\rho+\alpha)\|z^{\nu + 1} - z^{\nu}\|^{2}.
\end{aligned}
\label{eq:Ldiffzfinal}
\end{equation}
\normalsize
By substituting (\ref{eq:Ldiffxfinal}), (\ref{eq:Ldifflambda}) and (\ref{eq:Ldiffzfinal}) into (\ref{eq:Ldifference}), we obtain
\begin{equation}
\begin{aligned}
&L(x^{\nu+1}, z^{\nu+1}, \lambda^{\nu+1})-L(x^{\nu}, z^{\nu}, \lambda^{\nu})\\
& \leq\frac{\gamma M_{2}^{2} - \rho}{2 M_{2}^{2}} \sum_{k\in\mathcal{A}_{\nu}}\|x_{k}^{\nu+1}-x_{k}^{\nu}\|^2-(\rho+\alpha)\|z^{\nu + 1} - z^{\nu}\|^{2}\\&~~~~+\frac{1}{\rho}\sum_{k\in\mathcal{A}_{\nu}}\|\lambda_{k}^{\nu+1} - \lambda_{k}^{\nu}\|^{2} + \sum_{k\in\mathcal{A}_{\nu}}(\lambda_{k}^{\nu+1} - \lambda_{k}^{\nu})^{\top}(z_{k}^{\bar{\nu}_{k}+1}- z_{k}^{\nu})\\
&~~~~+\rho\sum_{k\in\mathcal{A}_{\nu}}({z}_{k}^{\bar{\nu}_{k}+1} - {z}_{k}^{\nu})^{\top}(A_{k}x_{k}^{\nu+1}-A_{k}x_{k}^{\nu}) \\
&\leq \frac{\gamma M_{2}^{2} - \rho}{2 M_{2}^{2}} \sum_{k\in\mathcal{A}_{\nu}}\|x_{k}^{\nu+1}-x_{k}^{\nu}\|^2-(\rho+\alpha)\|z^{\nu + 1} - z^{\nu}\|^{2}\\
&~~~~+(\frac{1}{\rho}+\frac{1}{2})\sum_{k\in\mathcal{A}_{\nu}}\|\lambda_{k}^{\nu+1}- \lambda_{k}^{\nu}\|^{2} + \frac{1}{2}\sum_{k\in\mathcal{A}_{\nu}}\|z_{k}^{\bar{\nu}_{k}+1}- z_{k}^{\nu}\|^{2} \\
&~~~~ + \rho M_{2}^{4} \sum_{k\in\mathcal{A}_{\nu}}\|{z}_{k}^{\bar{\nu}_{k}+1} - {z}_{k}^{\nu}\|^{2}
+\frac{\rho}{4 M_{2}^{4}}\sum_{k\in\mathcal{A}_{\nu}}\|A_{k}x_{k}^{\nu+1}-A_{k}x_{k}^{\nu}\|^{2}\\
& \leq (\frac{\gamma}{2} - \frac{\rho}{4 M_{2}^{2}}) \sum_{k\in\mathcal{A}_{\nu}}\|x_{k}^{\nu+1}-x_{k}^{\nu}\|^2-(\rho+\alpha)\|z^{\nu + 1} - z^{\nu}\|^{2}\\
&~~~~ +(\frac{1}{\rho}+\frac{1}{2})\sum_{k\in\mathcal{A}_{\nu}}\|\lambda_{k}^{\nu+1}- \lambda_{k}^{\nu}\|^{2} +\frac{2\rho M_{2}^{4}+1}{2} \sum_{k\in\mathcal{A}_{\nu}}\|{z}_{k}^{\bar{\nu}_{k}+1} - {z}_{k}^{\nu}\|^{2}.
\end{aligned}
\label{eq:Ldiffinter}
\end{equation}
\normalsize
The second inequality is obtained by applying Young's inequality; i.e., 
\small$ab \leq \frac{a^{2}}{2 \epsilon} + \frac{\epsilon b^{2}}{2}$\normalsize, with \small$\epsilon = 1~$\normalsize and \small$\epsilon = \frac{1}{2M_{2}^{4}}$\normalsize, to term \small$(\lambda_{k}^{\nu+1} - \lambda_{k}^{\nu})^{\top}(z_{k}^{\bar{\nu}_{k}+1}- z_{k}^{\nu})$\normalsize and \small$\rho({z}_{k}^{\bar{\nu}_{k}+1} - {z}_{k}^{\nu})^{\top}(A_{k}x_{k}^{\nu+1}-A_{k}x_{k}^{\nu})$\normalsize, respectively. The last inequality holds under Assumption \ref{assumption:Axmapping}.

We further bound \small$\| \lambda_{k}^{\nu+1} - \lambda_{k}^{\nu}\|^{2}$\normalsize as follows. From the optimality condition of (\ref{eq:xupdateglobal}) and (\ref{eq:lambdaupdateglobal}), for $k \in \mathcal{A}_{\nu}$, we have
\begin{equation}
\begin{aligned}
\bm{0} &= \partial F_{k}(x_{k}^{\nu+1})+A_{k}^{\top}\bigg(\lambda_{k}^\nu+\rho(A_{k}{x}_{k}^{\nu+1}-{z}_{k}^{\bar{\nu}_{k}+1})\bigg)\\
&=\partial F_{k}(x_{k}^{\nu+1})+A_{k}^{\top} \lambda_{k}^{\nu+1}.
\end{aligned}
\label{eq:xoptimality}
\end{equation}
\normalsize

For \small$\forall k \notin \mathcal{A}_{\nu}$\normalsize, let \small$\bar{\nu}_{k} < \nu ~$\normalsize denote the last iteration number for which region $k$ updated. Then, \small$\partial F_{k}(x_{k}^{\bar{\nu}_{k}})+A_{k}^{\top} \lambda_{k}^{\bar{\nu}_{k}+1}=\bm{0}$\normalsize. Since \small$x_{k}^{\bar{\nu}_{k}+1} = \cdots = x_{k}^{\nu} = x_{k}^{\nu+1}$\normalsize and \small$\lambda_{k}^{\bar{\nu}_{k}+1} = \cdots = \lambda_{k}^{\nu} = \lambda_{k}^{\nu+1}$\normalsize, we can conclude that
\begin{equation}
\partial F_{k}(x_{k}^{\nu+1})+A_{k}^{\top} \lambda_{k}^{\nu+1} = \bm{0},~~~~~ \forall k ~\text{and}~\forall \nu.
\label{eq:lambdaoptimality}
\end{equation}
\normalsize
By Assumption \ref{assumption:proxregular}, Assumption \ref{assumption:Ainvertible} and (\ref{eq:lambdaoptimality}), we can bound
\begin{equation}
\begin{aligned}
&\| \lambda_{k}^{\nu+1} - \lambda_{k}^{\nu}\|^{2} = \|(A_{k}A_{k}^{\top})^{-1}A_{k}(A_{k}^{\top}\lambda_{k}^{\nu+1}-A_{k}^{\top}\lambda_{k}^{\nu})\|^{2} \\
&= \|(A_{k}A_{k}^{\top})^{-1}A_{k}\big(\partial F_{k}(x_{k}^{\nu+1}) - \partial F_{k}(x_{k}^{\nu})\big)\|^{2}\\
&= \|B_{k}\big(\partial F_{k}(x_{k}^{\nu+1}) - \partial F_{k}(x_{k}^{\nu})\big)\|^{2}\\
& \leq \sigma_{\max}(B_{k}^{\top}B_{k}) \|\partial F_{k}(x_{k}^{\nu+1}) - \partial F_{k}(x_{k}^{\nu})\|^{2}\\
&\leq C M_{1}^2 \|x_{k}^{\nu+1}-x_{k}^{\nu}\|^2.
\end{aligned}
\label{eq:lambdabound}
\end{equation}
\normalsize
Substituting (\ref{eq:lambdabound}) into (\ref{eq:Ldiffinter}), we obtain Lemma \ref{lemma:Ldifference}.
\label{proof:Ldifference}
\end{proof}
\vspace{-0.1cm}
\section{Proof of Lemma~\ref{lemma:errorbound}}
\begin{proof}[Proof of Lemma~\ref{lemma:errorbound}]
\allowdisplaybreaks
\begin{equation}
\begin{aligned}
&\sum_{\phi = 1}^{\nu}\sum_{k\in\mathcal{A}_{\nu}}\|z_{k}^{\phi} - z_{k}^{\bar{\phi}_{k}+1}\|^{2} = \sum_{\phi = 1}^{\nu} \sum_{k\in\mathcal{A}_{\nu}}\bigg\|\sum_{\iota = \bar{\phi}_{k}+1}^{\phi - 1} (z_{k}^{\iota + 1}- z_{k}^{\iota})\bigg\|^{2}\\
& \leq \sum_{\phi = 1}^{\nu} \sum_{k\in\mathcal{A}_{\nu}}(\phi -\bar{\phi}_{k} -1) \sum_{\iota = \bar{\phi}_{k}+1}^{\phi - 1} \|z_{k}^{\iota + 1}- z_{k}^{\iota}\|^{2} \\
& \leq \sum_{\phi = 1}^{\nu} \sum_{k\in\mathcal{A}_{\nu}} (\omega -1) \sum_{\iota = \phi - \omega+1}^{\phi - 1} \|z_{k}^{\iota + 1}- z_{k}^{\iota}\|^{2}\\
& \leq (\omega-1)^{2} \sum_{\phi = 1}^{\nu} \sum_{k = 1}^{K}\|z_{k}^{\phi + 1}- z_{k}^{\phi}\|^{2}\\
& \leq (\omega-1)^{2} \sum_{\phi = 1}^{\nu} \|z^{\phi + 1}- z^{\phi}\|^{2},
\end{aligned}
\end{equation}
\normalsize
where for the second inequality, we apply (\ref{eq:iterdelaybound}); i.e., $\phi - \omega \leq \bar{\phi}_{k} < \phi$. The third inequality holds due to the observation that in the summation $\sum_{\phi = 1}^{\nu} \sum_{\iota = \phi - \omega+1}^{\phi - 1} \|z_{k}^{\iota + 1}- z_{k}^{\iota}\|^{2}$, each $\|z_{k}^{\iota + 1}- z_{k}^{\iota}\|$ with $\iota = 1,2, ..., \nu$ appears no more than $\omega - 1$ times. 
\label{proof:errorbound}
\end{proof}
\bibliographystyle{IEEEtran}
\bibliography{Partition,Partition2,CommunicationinSG,Asynchronous}

\begin{thebibliography}{10}
\providecommand{\url}[1]{#1}
\csname url@samestyle\endcsname
\providecommand{\newblock}{\relax}
\providecommand{\bibinfo}[2]{#2}
\providecommand{\BIBentrySTDinterwordspacing}{\spaceskip=0pt\relax}
\providecommand{\BIBentryALTinterwordstretchfactor}{4}
\providecommand{\BIBentryALTinterwordspacing}{\spaceskip=\fontdimen2\font plus
\BIBentryALTinterwordstretchfactor\fontdimen3\font minus
  \fontdimen4\font\relax}
\providecommand{\BIBforeignlanguage}[2]{{%
\expandafter\ifx\csname l@#1\endcsname\relax
\typeout{** WARNING: IEEEtran.bst: No hyphenation pattern has been}%
\typeout{** loaded for the language `#1'. Using the pattern for}%
\typeout{** the default language instead.}%
\else
\language=\csname l@#1\endcsname
\fi
#2}}
\providecommand{\BIBdecl}{\relax}
\BIBdecl

\bibitem{kar2014distributed}
S.~Kar, G.~Hug, J.~Mohammadi, and J.~M.~F. Moura, ``Distributed state
  estimation and energy management in smart grids: A consensus {+} innovations
  approach,'' \emph{IEEE Journal of Selected Topics in Signal Processing},
  vol.~8, no.~6, pp. 1022--1038, 2014.

\bibitem{conejo2006decomposition}
A.~J. Conejo, E.~Castillo, R.~Garc{\'\i}a-Bertrand, and R.~M{\'\i}nguez,
  \emph{Decomposition techniques in mathematical programming: engineering and
  science applications}.\hskip 1em plus 0.5em minus 0.4em\relax Springer
  Berlin, 2006.

\bibitem{bertsekas1989parallel}
D.~P. Bertsekas and J.~N. Tsitsiklis, \emph{Parallel and distributed
  computation: numerical methods}.\hskip 1em plus 0.5em minus 0.4em\relax
  Prentice hall Englewood Cliffs, NJ, 1989, vol.~23.

\bibitem{magnusson2016convergence}
S.~Magn{\'u}sson, P.~C. Weeraddana, M.~G. Rabbat, and C.~Fischione, ``On the
  convergence of alternating direction lagrangian methods for nonconvex
  structured optimization problems,'' \emph{IEEE Transactions on Control of
  Network Systems}, vol.~3, no.~3, pp. 296--309, 2016.

\bibitem{wang2015global}
Y.~Wang, W.~Yin, and J.~Zeng, ``Global convergence of admm in nonconvex
  nonsmooth optimization,'' \emph{arXiv preprint arXiv:1511.06324}, 2016.

\bibitem{boyd2011distributed}
S.~Boyd, N.~Parikh, E.~Chu, B.~Peleato, and J.~Eckstein, ``Distributed
  optimization and statistical learning via the alternating direction method of
  multipliers,'' \emph{Foundations and Trends{\textregistered} in Machine
  Learning}, vol.~3, no.~1, pp. 1--122, 2011.

\bibitem{dwork1988consensus}
C.~Dwork, N.~Lynch, and L.~Stockmeyer, ``Consensus in the presence of partial
  synchrony,'' \emph{Journal of the ACM (JACM)}, vol.~35, no.~2, pp. 288--323,
  1988.

\bibitem{peng2016arock}
Z.~Peng, Y.~Xu, M.~Yan, and W.~Yin, ``Arock: an algorithmic framework for
  asynchronous parallel coordinate updates,'' \emph{SIAM Journal on Scientific
  Computing}, vol.~38, no.~5, pp. A2851--A2879, 2016.

\bibitem{7285913}
I.~Aravena and A.~Papavasiliou, ``A distributed asynchronous algorithm for the
  two-stage stochastic unit commitment problem,'' in \emph{2015 IEEE Power
  Energy Society General Meeting}, July 2015, pp. 1--5.

\bibitem{abboud2014asynchronous}
A.~Abboud, R.~Couillet, M.~Debbah, and H.~Siguerdidjane, ``Asynchronous
  alternating direction method of multipliers applied to the direct-current
  optimal power flow problem,'' in \emph{IEEE International Conference on
  Acoustics, Speech and Signal Processing (ICASSP)}, 2014, pp. 7764--7768.

\bibitem{nguyen2016distributed}
H.~K. Nguyen, A.~Khodaei, and Z.~Han, ``Distributed algorithms for peak ramp
  minimization problem in smart grid,'' in \emph{IEEE International Conference
  on Smart Grid Communications (SmartGridComm)}, 2016, pp. 174--179.

\bibitem{chang2017scheduled}
C.-Y. Chang, J.~Cort{\'e}s, and S.~Mart{\'\i}nez, ``A scheduled-asynchronous
  distributed optimization algorithm for the optimal power flow problem,'' in
  \emph{American Control Conference (ACC), 2017}.\hskip 1em plus 0.5em minus
  0.4em\relax IEEE, 2017, pp. 3968--3973.

\bibitem{7423789}
T.~H. Chang, M.~Hong, W.~C. Liao, and X.~Wang, ``Asynchronous distributed admm
  for large-scale optimization -- part i: Algorithm and convergence analysis,''
  \emph{IEEE Transactions on Signal Processing}, vol.~64, no.~12, pp.
  3118--3130, June 2016.

\bibitem{kumar2017asynchronous}
S.~Kumar, R.~Jain, and K.~Rajawat, ``Asynchronous optimization over
  heterogeneous networks via consensus admm,'' \emph{IEEE Transactions on
  Signal and Information Processing over Networks}, vol.~3, no.~1, pp.
  114--129, 2017.

\bibitem{cannelli2016asynchronous}
L.~Cannelli, F.~Facchinei, V.~Kungurtsev, and G.~Scutari, ``Asynchronous
  parallel algorithms for nonconvex big-data optimization-part i: Model and
  convergence,'' \emph{Submitted to Mathematical Programming}, 2016.

\bibitem{zhang2014asynchronous}
R.~Zhang and J.~Kwok, ``Asynchronous distributed admm for consensus
  optimization,'' in \emph{International Conference on Machine Learning}, 2014,
  pp. 1701--1709.

\bibitem{erseghe2015distributed}
T.~Erseghe, ``A distributed approach to the \uppercase{OPF} problem,''
  \emph{EURASIP Journal on Advances in Signal Processing}, vol. 2015, no.~1,
  pp. 1--13, 2015.

\bibitem{guo2016acase}
J.~Guo, G.~Hug, and O.~K. Tonguz, ``A case for nonconvex distributed
  optimization in large-scale power systems,'' \emph{IEEE Transactions on Power
  Systems}, vol.~32, no.~5, pp. 3842--3851, 2017.

\bibitem{chiang2017feasible}
H.-D. Chiang and C.-Y. Jiang, ``Feasible region of optimal power flow:
  Characterization and applications,'' \emph{IEEE Transactions on Power
  Systems}, 2017.

\bibitem{guointelligent}
J.~Guo, G.~Hug, and O.~K. Tonguz, ``Intelligent partitioning in distributed
  optimization of electric power systems,'' \emph{IEEE Transactions on Smart
  Grid}, vol.~7, no.~3, pp. 1249--1258, 2016.

\end{thebibliography}

\end{document}